\documentclass{amsart}
\usepackage{amsmath,amsthm,latexsym,amssymb}
\usepackage{eucal}
\usepackage{mathrsfs}
\usepackage[all]{xy}
\allowdisplaybreaks

\numberwithin{equation}{section}

\newtheorem{thm}{Theorem}[section]
\newtheorem*{thm*}{Theorem}
\newtheorem{lem}[thm]{Lemma}
\newtheorem{prop}[thm]{Proposition}
\newtheorem{cor}[thm]{Corollary}

\theoremstyle{definition}
\newtheorem{defn}[thm]{Definition}
\newtheorem{rmk}[thm]{Remark}
\newtheorem*{rmk*}{Remark}
\newtheorem{ex}[thm]{Example}

\theoremstyle{remark}

\newtheorem{notn}[thm]{Conventions}

\newdir{ >}{{}*!/-10pt/@{>}}
\newcommand\fib{\ar @{->>} [r]} 
\newcommand\cof{\ar @{ >->}[r]}

\newcommand \Om{\Omega}
\newcommand \n{\nabla}
\newcommand \del{\partial}
\newcommand \vp{\varphi}

\newcommand{\sdr}[2]{\overset{#1}{\underset{#2}{\rightleftharpoons}}}

\newcommand \sn[1]{(-1)\sp{#1}}
\renewcommand\H{\operatorname{H}}
\newcommand\op{\mathcal}
\newcommand\cat{\mathbf}
\newcommand{\ob}{\operatorname{Ob}}

\newcommand{\hoch}{\mathscr H}

\renewcommand{\Bar}{\mathscr B}
\newcommand{\cohoch}{\widehat{\mathscr H}}

\newcommand{\si}{s^{-1}}

\newcommand\E{\mathsf E}

\begin{document}
\title [CoHochschild homology]{CoHochschild homology of chain coalgebras}
\author{Kathryn Hess}
\author{Paul-Eug\`ene Parent}
\author{Jonathan Scott}

\address{Institut de g\'eom\'etrie, alg\`ebre et topologie (IGAT) \\
    \'Ecole Polytechnique F\'ed\'erale de Lausanne \\
    CH-1015 Lausanne \\
    Switzerland}
    \email{kathryn.hess@epfl.ch}
\address{Department of Mathematics and Statistics \\
    University of Ottawa \\
    585 King Edward Avenue \\
    Ottawa, ON \\
    K1N 6N5 Canada}
    \email{pparent@uottawa.ca}
\address{Department of Mathematics and Statistics \\
    University of Ottawa \\
    585 King Edward Avenue \\
    Ottawa, ON \\
    K1N 6N5 Canada}
    \email{jscott@uottawa.ca}
\date { \today}

  \begin{abstract} Generalizing work of Doi and of Idrissi, we define a coHochschild homology theory for chain coalgebras over any commutative ring and prove its naturality with respect to morphisms of chain coalgebras up to strong homotopy.  As a consequence we obtain that if the comultiplication of a chain coalgebra $C$ is itself a morphism of chain coalgebras up to strong homotopy, then the coHochschild complex $\cohoch (C)$  admits a natural comultiplicative structure.  In particular, if $K$ is a reduced simplicial set and $C_{*}K$ is its normalized chain complex, then $\cohoch (C_{*}K)$ is naturally a homotopy-coassociative chain coalgebra.  We provide a simple, explicit formula for the comultiplication on $\cohoch (C_{*}K)$ when $K$ is a simplicial suspension.
  
The coHochschild complex construction is topologically relevant. Given two simplicial maps $g,h:K\to L$, where $K$ and $L$ are reduced, the homology of the coHochschild complex of $C_{*}L$ with coefficients in $C_{*}K$ is isomorphic to the   homology of the homotopy coincidence space of the geometric realizations of $g$ and $h$, and this isomorphism respects comultiplicative structure. In particular,  there a isomorphism, respecting comultiplicative structure, from the homology of $\cohoch(C_{*}K)$ to $H_{*}\op L|K|$, the homology of the free loops on the geometric realization of $K$.
 \end{abstract}
 
 \maketitle

\tableofcontents
 

\section* {Introduction}\label{sec:introduction}

Hochschild homology is a well-known and very useful homology theory for algebras, which has considerable relevance in topology as well.  In particular, for any based topological space $X$, the Hochschild homology of $S_{*}(\Om X)$, the singular chains on the space of based loops on $X$, is isomorphic to the singular homology of the space $\op L X$  of free loops on $X$ (cf, e.g., \cite{jones, loday}).

In \cite{doi} Doi developed a homology theory for coalgebras over a field that is analogous to the Hochschild homology of algebras.  In this article, we offer an alternate approach to Doi's homology theory, which allows us to extend his theory, which we call \emph{coHochschild homology}, easily to chain coalgebras over any commutative ring. In particular, we describe
 the coHochschild complex $\cohoch (N,C)$ of a chain coalgebra $C$ with coefficients in a bicomodule $N$ as a twisted extension of the cobar construction on $C$.

We prove that the coHochschild complex $\cohoch (C)$ is natural with respect to morphisms of chain coalgebras up to strong homotopy (Theorem \ref{thm:ext-nat-cohoch}). Together with the fact that the coHochschild functor is comonoidal (Theorem \ref{thm:milgram-cohoch}), this extended naturality enables us to prove that if the comultiplication on $C$ is itself a morphism of coalgebras up to strong homotopy, then $\cohoch(C)$ admits a natural comultiplication.  Moreover, we determine conditions under which this natural comultiplication is coassociative, either strictly or up to chain homotopy (Theorem \ref{thm:cohoch-comult}).  We also establish a more general version of this result, for $\cohoch(C,C')$, where $C$ is a chain coalgebra seen as a $C'$-bicomodule via two chain coalgebra maps $f,g:C\to C'$ (Theorem \ref{thm:rel-cohoch-comult}). 

Let $K$ be a reduced simplicial set, and let $C_{*}K$ denote the normalized chain complex on $K$. It follows from the purely algebraic results cited above and from our earlier work (\cite {hpst},\cite {hps2}) that $\cohoch (C_{*}K)$ admits a canonical comultiplication, which is coassociative up to chain homotopy. We provide a simple, explicit formula for this comultiplication when $K$ is a simplicial suspension (Example \ref{ex:comult-susp}).  More generally, there is a simple formula for the comultiplication in $\cohoch(C_{*}K,C_{*}L)$, where $K$ is a simplicial suspension, $L$ is a reduced simplicial set, and  the $C_{*}L$-bicomodule structure on $C_{*}K$ is determined by two simplicial maps $g,h:K\to L$ (Example \ref{ex:comult-coince}).

We illustrate the topological utility of the coHochschild complex and its comultiplicative structure, when we prove the theorem below, concerning a certain homotopy-invariant version of the coincidence space of two continuous maps $g,h:X\to Y$:
 $$\op E_{g,h}=\big\{(x,\ell)\in X\times Y^{I}\mid \ell (0)=g(x), \ell (1)=h(x)\big\}.$$
Note that the free loop space on $X$ is just $\op E_{Id_{X},Id_{X}}$.  More generally, let $Y$ be a manifold, and let $U$ and $V$ be submanifolds of $Y$.  Let $g:U\times V\to Y$ be given by projection onto the first coordinate, while $h:U\times V\to Y$ is given by projection onto the second coordinate.  The homotopy coincidence space $\op E_{g,h}$ is then exactly the space of open strings in $Y$ starting in $U$ and ending in $V$. 

\begin{thm*}[Theorem \ref{thm:model-htpy-coinc}] If $g,h:K\to L$ are simplicial maps, where $L$ is a reduced simplicial set, then there is a quasi-isomorphism of chain complexes
$$\cohoch (C_{*}K,C_{*}L)\xrightarrow\simeq S_{*}\op E_{|g|,|h|}$$
that is comultiplicative up to chain homotopy.
\end{thm*}

Here $|g|$ and $|h|$ denote the geometric realizations of the simplicial maps.

If $K$ and $L$ have only finitely many nondegenerate simplices, then $\cohoch (C_{*}K,C_{*}L)$ is a finitely generated module over a finitely generated, free algebra, endowed with a relatively simple differential and an explicitly defined comultiplication.  It is thus realistic to expect to be able to make explicit homology computations, including comultiplicative structure, with this model for the homotopy coincidence space.

In the appendix to this article we show that for any pair of reduced simplicial coalgebras $M_{\bullet}$ and $M'_{\bullet}$ over the base ring $R$, the natural chain equivalence
$$f:\op A_{N}( M_{\bullet}\boxtimes M'_{\bullet})\to \op A_{N}(M_{\bullet})\otimes \op A_{N}(M'_{\bullet})$$
is strongly homotopy comultiplicative, where $\boxtimes$ denotes levelwise tensor product over $R$ and $\op A_{N}$ denotes the normalized chain complex functor (Theorem \ref{thm:aw-dcsh}).  As a consequence we obtain a generalization of the main results of \cite{hpst}, establishing that the cobar construction on the normalized chain complex of any reduced simplicial set admits a comultiplication that is at least homotopy-coassociative.

\begin{rmk*} The coHochschild complex  of a chain coalgebra plays a very important role in \cite{hess-rognes}, where it is the essential building block in the construction  of a chain complex model for the spectrum homology of topological cyclic homology of a topological space.  The power maps on the coHochschild complex of $C_{*}K$, as constructed in \cite {hess-rognes:power}, which are algebraic models for the topological power maps on the free loop space, are the key elements of this construction.

We note further that  Theorem \ref{thm:model-htpy-coinc} is crucial to the proof in  \cite{hess-rognes:power} that the algebraic power map on $\cohoch (C_{*}K)$ is indeed a model for the topological power map.
 \end{rmk*}

\begin{rmk*}  Let $K$ be any reduced simplicial set. Dualizing the coHochschild complex $\cohoch (C_{*}K)$ and its homotopy-coassociative comultiplication, one obtains a homotopy-associative multiplication on the Hochschild cochain complex  for $C^*K$ with coefficients in $C_{*}K$.  The nature of this dualized multiplication is entirely different from that of the well-known multiplication on the Hochschild cochain complex  for $C^*K$ with coefficients in $C^{*}K$, which is of purely algebraic origin.
\end{rmk*}

\begin{rmk*} In \cite{idrissi} Idrissi sketched a proof of the existence of a homotopy-coassociative comultiplication on the coHochschild complex of any chain coalgebra over a field, with comultiplication that is a morphism of coalgebras up to strong homotopy.  The proof of one crucial lemma (Lemma 1) does not seem to be complete, though. As formulated,  the author's proof requires extended naturality of the coHochschild construction,  to prove the existence of a map that is one of the factors of the purported comultiplication.  There is no proof of extended naturality in \cite {idrissi}, however.
\end{rmk*} 

The authors would like to thank the referee for his comments, which enabled us to sharpen the focus of this article.

\subsection*{Notation and conventions}
\begin{itemize}

\item Throughout this paper we are working over a principal ideal domain $R$.  We denote the category of graded $R$-modules by $\cat{grMod}_R$, the category of chain complexes over $R$ by $\cat{Ch}_R$, the category of augmented chain algebras over $R$ by $\cat {Alg}_{R}$ and the category of coaugmented, connected chain coalgebras by $\cat{Coalg}_{R}$.  The underlying graded modules of all chain (co)algebras are assumed to be $R$-free.  
\item The degree of an element $v$ of a graded module $V$ is denoted $|v|$. 
\item Throughout this article we apply the Koszul sign convention for commuting elements  of a graded module or for commuting a morphism of graded modules past an element of the source module.  For example,  if $V$ and $W$ are graded algebras and $v\otimes w, v'\otimes w'\in V\otimes W$, then $$(v\otimes w)\cdot (v'\otimes w')=(-1)^{|w|\cdot |v'|}vv'\otimes ww'.$$ Furthermore, if $f:V\to V'$ and $g:W\to W'$ are morphisms of graded modules, then for all $v\otimes w\in V\otimes W$, 
$$(f\otimes g)(v\otimes w)=(-1)^{|g|\cdot |v|} f(v)\otimes g(w).$$
\item The \emph {suspension} endofunctor $s$ on the category of graded modules is defined on objects $V=\bigoplus _{i\in \mathbb Z} V_ i$ by
$(sV)_ i \cong V_ {i-1}$.  Given a homogeneous element $v$ in
$V$, we write $sv$ for the corresponding element of $sV$. The suspension $s$ admits an obvious inverse, which we denote $\si$.
\item Given chain complexes $(V,d)$ and $(W,d)$, the notation
$f:(V,d)\xrightarrow{\simeq}(W,d)$ indicates that $f$ induces an isomorphism in homology. 
In this case we refer to $f$ as a \emph {quasi-isomorphism}.
\item Let $f,g:A\to A'$ be morphisms of chain algebras.  A \emph{derivation homotopy} from $f$ to $g$ consists of a chain homotopy $H:A\to A'$ from $f$ to $g$ such that $H(ab)=H(a)f(b)+(-1)^{|a|}g(a) H(b)$ for all $a,b\in A$.
\item Let $T$ denote the endofunctor on the category of free graded $R$-modules given by
$$TV=\oplus _{n\geq 0}V^{\otimes n},$$
where $V^{\otimes 0}=R$.  An element of the summand $V^{\otimes n}$ of $TV$ is denoted $v_{1}|\cdots |v_{n}$, where $v_{i}\in V$ for all $i$. 
\item Let $N:TV\to TV$ denote the \emph{norm} operator given by
$$N(v_{1}|\cdots |v_{n})=\sum _{1\leq j\leq n}\pm v_{j}|\cdots |v_{n}|v_{1}|\cdots |v_{j-1},$$
the signed sum of cyclic permutations, where the sign is determined by the Koszul rule and by the sign of the permutation.  
\item If $A$ is an augmented chain algebra, then $\overline A$ denotes its augmentation ideal.  Similarly, the coaugmentation coideal of a coaugmented chain coalgebra $C$ is denoted $\overline C$.
\item The normalized chains functor from simplicial sets to chain complexes is denoted $C_{*}$, while the singular simplices functor from topological spaces to simplicial sets is denote $\mathcal S_{\bullet}$.  Their composite, $C_{*}\circ \mathcal S_{\bullet}$, is denoted $S_{*}$.  The left adjoint to $S_{\bullet}$, i.e.,  geometric realization, is denoted $|-|$.
\end{itemize}


\section{The cobar and coHochschild complexes of a chain coalgebra}\label{sec:cohoch}

In this section we introduce the coHochschild complex of a chain coalgebra over a principal ideal domain, generalizing the definitions in \cite{doi} and in \cite{idrissi}.  The prefix ``co'' in the name of this complex is justified by the fact that there is an underlying cosimplicial object in the category of chain complexes, which we do not define explicitly here.  

We begin by recalling the classical bar and Hochschild complexes of a chain algebra $A$.
Though these constructions are well known, we consider it worthwhile to present them again briefly here, for two reasons.  First, our presentation, while fairly standard from the perspective of a topologist, is rather different from that with which algebraists are familiar.  Second, the dual constructions for chain coalgebras are easier to understand when compared directly with the known constructions for chain algebras.

All signs in the formulas below follow from the Koszul rule. It is a matter of straightforward calculation in each case to show that the differential squares to zero.

\subsection{The bar and Hochschild complexes}\label{sec:prelims}

Let $\Bar$ denote the \emph{bar construction} functor from $\cat{Alg}_{R}$ to $\cat {Ch}_{R}$, defined by 
$$\Bar A=\left(T (s\overline A), d_{\Bar}\right)$$
where, if $d$ is the differential on $A$, then
\begin{align*}
d_{\Bar}(sa_{1}|\cdots|sa_{n})=&\sum _{1\leq j\leq n}\pm sa_{1}|\cdots |s(da_{j})|\cdots |sa_{n}\\ 
&+\sum _{1\leq j<n}\pm sa_{1}|...|s(a_{j}a_{j+1})|\cdots |sa_{n}.
\end{align*}

Observe that the graded $R$-module underlying $\Bar A$ is naturally a cofree coassociative coalgebra, with comultiplication given by splitting of words. The differential $d_{\Bar}$ is a coderivation with respect to this splitting comultiplication, so that $\Bar A$ is itself a chain coalgebra.  If $C$ is a conilpotent chain coalgebra, then any chain coalgebra map $\gamma:C\to \Bar A$ is determined by its projection to the coalgebra cogenerators $s\overline A$, denoted $\gamma_{1}$.

Let $\cat {BiMod}$ denote the category in which the objects are pairs $(A,M)$, where $A$ is a connected, augmented chain $R$-algebra and $M$  a chain $A$-bimodule endowed with an augmentation $M\to R$ that is a morphism of $A$-bimodules. A morphism from $(A,M)$ to $(A',M')$ consists of a pair $(f,g)$, where $f:A\to A'$ is a morphism of chain algebras and $g:M\to M'$ is a morphism of aumented, chain $A$-bimodules with respect to the $A$-bimodule structure on $M'$ induced by $f$. 

As a lift of the bar construction, let $\hoch(-,-)$ denote the \emph{Hochschild complex} functor from $\cat{BiMod}$ to $\cat {Ch}_{R}$, defined by 
$$\hoch (A,M)= \left(T (s\overline A)\otimes M, d_{\hoch}\right)$$
where, if $d$ denotes the differentials on $A$ and on $M$, then
\begin{align*}
d_{\hoch}(sa_{1}|\cdots|sa_{n}\otimes x)=&d_{\Bar}(sa_{1}|\cdots|sa_{n})\otimes x \;\pm sa_{1}|\cdots|sa_{n}\otimes dx\\
&+sa_{1}|\cdots|sa_{n-1}\otimes a_{n}\cdot x\; \pm sa_{2}|\cdots |sa_{n}\otimes x\cdot a_{1},
\end{align*}
for all $sa_{1}|\cdots|sa_{n}\otimes x\in T (s\overline A)\otimes M$, where $\cdot$ denotes both the right and the left actions of $A$ on $M$.

For every object $(A,M)$ in $\cat {BiMod}$,  there  is clearly a twisted extension of chain complexes
\begin{equation}\label{eqn:hoch-ext}
\xymatrix@1{M\cof &\hoch (A,M)\fib &\Bar A.}
\end{equation}

When $A$ is considered as a bimodule over itself, where the bimodule structure is given by multiplication in $A$,  we write
$$\hoch (A):=\hoch (A,A).$$
We consider $\hoch(-)$ as a functor from $\cat {Alg}_{R}$ to $\cat {Ch}_{R}$.

\begin{rmk} Let $X$ be a based topological space.  It is well known (e.g., \cite{jones}) that the homology of $\hoch(S_{*}\Om X)$ is isomorphic to the homology of the space $\mathcal LX$ of free (i.e., unbased) loops on $X$.  In section \ref{sec:fls} we expand upon and generalize this result in the context of the dual constructions for chain coalgebras.
\end{rmk}

\subsection{The cobar construction}

Using fairly standard notation, let $\Om$ denote the \emph{cobar construction} functor from $\cat {Coalg}_{R}$ to $\cat {Ch}_{R}$, defined by 
$$\Om C= \left(T (\si \overline C), d_{\Om}\right)$$
where, if $d$ denotes the differential on $C$, then
\begin{align*}
d_{\Om}(\si c_{1}|\cdots|\si c_{n})=&\sum _{1\leq j\leq n}\pm \si c_{1}|\cdots |\si (dc_{j})|\cdots |\si c_{n}\\ 
&+\sum _{1\leq j\leq n}\pm \si c_{1}|...|\si c_{ji}|\si c_{j}{}^{i}|\cdots |\si c_{n},
\end{align*}
with signs determined by the Koszul rule, where the reduced comultiplication applied to $c_{j}$ is $c_{ji}\otimes c_{j}{}^{i}$ (using Einstein implicit summation notation).  

Observe that the graded $R$-module underlying $\Om C$ is naturally a free associative algebra, with multiplication given by concatenation. The differential $d_{\Om }$ is a derivation with respect to this concatenation product, so that $\Om C$ is itself a chain algebra.  Any chain algebra map $\alpha:\Om C\to A$ is determined by its restriction to the algebra generators $\si \overline C$. 

The cobar construction functor is comonoidal, i.e., there is a natural transformation $\Om (-\otimes -)\to \Om (-)\otimes \Om (-)$ that is appropriately coassociative and counital.  Milgram defined this natural transformation in \cite{milgram} for simply connected chain coalgebras and showed that it was a quasi-isomorphism.  In \cite{hps2} the authors extended the definition to all coaugmented chain coalgebras and showed that it was  actually a natural chain homotopy equivalence.

\begin{thm}\label{thm:milgram-sdr}
Let $C$ and $C'$ be coaugmented chain coalgebras.
There is a chain homotopy equivalence 
$$q:\Om (C\otimes C')\xrightarrow \simeq \Om C\otimes \Om C',$$
specified by 
$q\big(\si (c\otimes 1)\big)=\si c\otimes 1$,  $q\big(\si (1\otimes c')\big)=1\otimes \si c'$ and $q\big(\si (c\otimes c')\big)=0$ if $c\in \overline C$ and $c'\in \overline C'$.
\end{thm} 

Seen as functors from coalgebras to algebras and vice-versa, the cobar and bar constructions form an adjoint pair $\Om \dashv \Bar$.  Let $\eta: Id\to \Bar \Om$ denote the unit of this adjunction.  It is well known that for all coaugmented chain coalgebras $C$, the unit map
\begin{equation}\label{eqn:unit-barcobar}
\eta_{C}:C\xrightarrow\simeq\Bar\Om C
\end{equation}
is a quasi-isomorphism of chain coalgebras.  Furthermore, $\eta_{C}$ admits a natural retraction (i.e., a left inverse)
\begin {equation}\label{eqn:retract-barcobar}
\rho_{C}:\Bar\Om C\xrightarrow \simeq C
\end{equation}
that is a morphism of chain complexes, but a morphism of chain coalgebras only up to strong homotopy.

\subsection{The coHochschild complex}
Let $\cat {BiComod}$ denote the category in which the objects are pairs $(N,C)$, where $C$ is a coaumented chain $R$-coalgebra and $N$ is a chain $C$-bicomodule such that $N$ admits a coaugmentation $R\to N$ that is a morphism of $C$-bicomodules.  A morphism from $(C,N)$ to $(C',N')$ consists of a pair $(g,f)$, where $f:C\to C'$ is a morphism of chain coalgebras and $g:N\to N'$ is a morphism of chain $C'$-bicomodules with respect to the $C'$-bicomodule structure on $N$ induced by $f$. 

As an extension of the cobar construction, let $\cohoch(-,-)$ denote the \emph{coHochschild complex} functor from $\cat {BiComod}$ to $\cat {Ch}_{R}$, defined as follows.  
Let $C$ be a connected, coaugmented chain coalgebra, and let $N$ be  a $C$-bicomodule, with coactions $\lambda : N\to C\otimes N$ and $\rho:N\to N\otimes C$.  Applying the Einstein implicit-summation convention, write $\lambda (x)=e_{i}\otimes x^{i}$ and $ \rho(x)=x_{j}\otimes e^j$.  We then let
$$\cohoch(N,C)= \left(N\otimes T (\si \overline C), d_{\cohoch}\right)$$
where
\begin{align*}
d_{\cohoch}(x\otimes \si c_{1}|\cdots|\si c_{n})=&dx\otimes \si c_{1}|\cdots|\si c_{n}\;\pm x\otimes d_{\Om}(\si c_{1}|\cdots|\si c_{n})\\
&\pm x_{j}\otimes \si e^{j}|\si c_{1}|\cdots|\si c_{n}\\
& \pm x^{i}\otimes \si c_{1}|\cdots|\si c_{n}|\si e_{i},
\end{align*}
with the convention that applying $\si$ to an element of degree 0 gives $0$. The signs follow from the Koszul rule, as usual.

For every object $(N,C)$ in $\cat {BiComod}$, there is clearly a twisted extension of chain complexes
\begin{equation}\label{eqn:cohoch-ext}
\xymatrix@1{\Om C\cof &\cohoch (N,C)\fib &N.}
\end{equation}

The coHochschild complex functor is comonoidal with respect to the obvious monoidal structure on the category $\cat {BiComod}$, via a natural chain homotopy equivalence
$$\cohoch(-\otimes -, -\otimes -) \to \cohoch (-,-)\otimes \cohoch (-,-)$$
that extends Milgram's natural transformation for the cobar construction (Theorem \ref{thm:milgram-sdr}).

\begin{thm}\label{thm:milgram-cohoch} For all $(N,C)$ and $(N', C')$, objects in $\cat{BiComod}$, there is a natural chain homotopy equivalence
$$\hat q: \cohoch(N\otimes N', C\otimes C') \to \cohoch (N,C)\otimes \cohoch (N',C')$$
such that 
$$\xymatrix{\Om (C\otimes C')\ar[d] \ar[r]^q_{\simeq}&\Om C\otimes \Om C'\ar[d]\\
\cohoch (N\otimes N',C\otimes C') \ar[d]\ar[r]^{\hat q}_{\simeq }&\cohoch (N,C)\otimes\cohoch (N',C')\ar[d]\\
N\otimes N'\ar@{=}[r]&N\otimes N'}$$
commutes, where the vertical arrows are the natural inclusions and projections. 
\end{thm}

\begin{proof}The proof in \cite[Appendix A]{hps2} can easily be generalized to this situation.  The morphism $\hat q$ of $\Om (C\otimes C')$-modules is specified by 
$$\hat q (x\otimes x' \otimes 1)= (x\otimes 1)\otimes (x'\otimes 1)$$ 
for all $x\in N$, $x'\in N'$, which clearly gives rise to a differential map when extended as a map from a free right $\Om (C\otimes C')$-module. There is a section of $\hat q$
$$\hat \sigma: \cohoch (N,C)\otimes \cohoch (N',C')\to \cohoch(N\otimes N', C\otimes C')$$
given by the composite
$$\xymatrix{(N\otimes \Om C)\otimes (N'\otimes \Om C')\ar[r]&\big(N\otimes \Om (C\otimes C')\big)\otimes \big(N'\otimes \Om (C\otimes C')\big)\ar [d]^\cong\\
&(N\otimes N')\otimes \Om (C\otimes C')\otimes \Om (C\otimes C')\ar [d]^{Id_{N\otimes N'}\otimes \mu}\\
&(N\otimes N')\otimes \Om (C\otimes C'),}$$
where the first arrow is the obvious inclusion, and $\mu$ is the multiplication map on $\Om (C\otimes C')$.  It is easy to see that $\hat q\hat\sigma$ is the identity.  To complete the proof, one defines a chain homotopy $\hat h$ on $\cohoch(N\otimes N', C\otimes C')$ from $\hat\sigma\hat q$ to the identity, extending the homotopy $h$ from section A.2 of \cite{hps2}.
\end{proof}

\begin{ex}\label{ex:cohoch}  The following special cases of the coHochschild complex are worthy of note.
\begin {enumerate}
\item Considering the ground ring $R$ as a trivial $C$-bicomodule, for $C$ any chain coalgebra, we obtain that $\cohoch (R,C)=\Om C$.
\item When $C$ is considered as a comodule over itself, where the bicomodule structure is given by the comultiplication on $C$ on both sides, we write 
$$\cohoch (C):= \cohoch (C,C).$$
We consider $\cohoch(-)$ as a functor from $\cat {Coalg}_{R}$ to $\cat {Ch}_{R}$.  We show in section \ref{sec:comult-cohoch} that $\cohoch(-)$ is actually natural with respect to a much larger class of morphisms (Theorem \ref{thm:ext-nat-cohoch}).
\item Any coaugmented chain coalgebra $C$ can be considered as a bicomodule over itself, where the left coaction is trivial, i.e., equal to the composite
$$C\xrightarrow \cong R\otimes C \xrightarrow {\eta\otimes Id_{C}} C\otimes C.$$
If ${}_{\eta}C_{\Delta}$ denotes $C$ endowed with this $C$-bicomodule structure, then  $\cohoch ({}_{\eta}C_{\Delta}, C)$ is the usual acyclic cobar construction on $C$.
\item More generally, if $N$ is a right $C$-comodule, with right coaction $\rho$, then it can be considered as a $C$-bicomodule with trivial left $C$-coaction 
$$N\cong R\otimes N \xrightarrow {\eta\otimes Id_{N}}C\otimes N.$$
If ${}_{\eta}N_{\rho}$ denotes $N$ endowed with this $C$-bicomodule structure, then  $\cohoch ({}_{\eta}N_{\rho}, C)$ is the usual one-sided cobar construction on $C$ with coefficients in $N$.
\item If $C$ is a chain coalgebra with comultiplication $\Delta$, then $C\otimes C$ is naturally a $C$-bicomodule, where the right $C$-coaction is $Id_{C}\otimes \Delta$ and the left $C$-coaction is $\Delta\otimes Id_{C}$.  The coHochschild complex $\cohoch (C\otimes C, C)$ is isomorphic to the two-sided cobar construction on $C$.
\item More generally, if $M$ is a left $C$-comodule with coaction $\lambda$ and $N$ is a right $C$-comodule with coaction $\rho$, then $M\otimes N$ is naturally a $C$-bicomodule, with left coaction $\lambda\otimes Id_{N}$ and right coaction $Id_{M}\otimes \rho$.  The coHochschild complex $\cohoch (M\otimes N, C)$ is isomorphic to the two-sided cobar construction on $C$ with coefficients in $N$ on the left and in $M$ on the right.
\end{enumerate}
\end{ex}

\begin{rmk}  There is a natural and straightforward extension of the coHochschild complex of a chain coalgebra to a \emph{cocyclic complex}, analogous to the extension of the Hochschild complex of a chain algebra to the cyclic complex.
\end{rmk} 


\section{Comultiplicative structure on the coHochschild complex}\label{sec:comult-cohoch}

In this section we study comultiplicative structure on the coHochschild complex, determining, in particular, under what conditions such structure naturally exists.

We begin by establishing an ``extended naturality'' result for the coHochschild complex. It is clear from the definition of the complex $\cohoch (C)$ that it is natural with respect to morphisms of coalgebras.  We show below that it is in fact natural with respect to the much larger class of coalgebra morphisms up to strong homotopy, first defined by Gugenheim and Munkholm in \cite{gugenheim-munkholm}.

\begin{defn}\cite{gugenheim-munkholm}  
Given $C,C'\in \ob \cat {Coalg}_{R}$, a chain map $f:C\to C'$ is called a \emph{DCSH map} or a \emph{morphism of chain coalgebras up to strong homotopy} if there is a chain algebra map $\omega: \Om C\to \Om C'$ such that $\omega (\si c)= \si f(c) + \text {higher-order terms}$ for all $c\in C$.  The chain algebra map $\omega$ is said to \emph{realize the strong homotopy structure} of the DCSH map $f$.
\end{defn}

\begin{rmk} \label{rmk:unfold-dcsh} A chain algebra map $\omega:\Om C\to \Om C'$ is determined by a set of $R$-linear maps
$$\{\omega_{k}:\overline C \to (\overline C')^{\otimes k}\mid k\geq 1\},$$
where
\begin {enumerate}
\item $\omega _{1}=f$;
\item $\omega _{k}$ is homogeneous of degree $k-1$ for all $k$; and
\item \begin{align*}
\omega _{k}d_{C}+(-1)^kd_{(C')^{\otimes k}}\omega _{k}=&\sum _{i+j=k}(\omega _{i}\otimes \omega _{j})\overline \Delta_{C} \\
&-\sum _{i+j=k-2}(-1)^{i}(Id_{C'}^{\otimes i}\otimes\Delta_{C'}\otimes Id_{C'}^{\otimes j})\omega _{k-1} 
\end{align*}
for all $k$.
\end{enumerate} 
The relation between $\omega $ and the family $\{\omega_{k}\}_{k}$ is that 
$$\omega (\si e)=\sum _{k\geq 1}(\si)^{\otimes k}\omega _{k}(e),$$
for all $e\in C$.
\end{rmk}

The extended naturality of the coHochschild complex construction with respect to DCSH maps, stated precisely in the next theorem, is the key to obtaining comultiplicative structure on the coHochschild complex of certain coalgebras.

\begin {thm}\label{thm:ext-nat-cohoch}  A DCSH map $f:C\to C'$ with a fixed choice of chain algebra map $\omega:\Om C\to \Om C'$ realizing its strong homotopy structure naturally induces a chain map
$$\widehat \omega:\cohoch(C) \to \cohoch(C')$$
such that 
$$\xymatrix{\Om C\ar[d] \ar[r]^\omega&\Om C'\ar[d]\\
\cohoch (C) \ar[d]\ar[r]^{\widehat \omega}&\cohoch (C')\ar[d]\\
C\ar [r]^f&C'}$$
commutes, where the vertical arrows are the natural inclusions and projections. 
\end {thm}

\begin{proof} Let $\{\omega_{k}:\overline C \to (\overline C')^{\otimes k}\mid k\geq 1\}$ denote the family of $R$-linear maps associated to $\omega:\Om C\to \Om C'$, as in Remark \ref{rmk:unfold-dcsh}.

Given $e\in C$, we use Sweedler-type notation and write
$$\omega_{k}(e)=e'_{(k,1)}\otimes \cdots \otimes e'_{(k,k)},$$
for all $k$, suppressing the summation.
For $e\otimes w\in C\otimes \Om C$, set
$$\widehat \omega(e\otimes w)=\sum _{k\geq 1\atop 1\leq i\leq k}\pm e'_{(k,i)}\otimes \si e'_{(k,i+1)}|\cdots |\si e'_{(k,k)}\cdot \omega (w)\cdot \si e'_{(k,1)}|\cdots|\si e'_{(k,i-1)},$$
where the sign is given by the Koszul rule.

To show that $\widehat \omega$ is a chain map, we proceed as follows. Neglecting terms arising from the internal differential on $C$,  the summands in $\widehat \omega d_{\cohoch}(e\otimes w)$ are in one-to-one correspondence with the summands of
$$N\Big(\sum_{n\geq 1}\sum _{k+l=n}(\omega_{k}\otimes \omega _{l})\Delta(e)\Big),$$
while the summands in $d_{\cohoch}\widehat \omega (e\otimes w)$ are in one-to-one correspondence with the summands of
$$N\Big(\sum_{n\geq 1}\sum _{k+l=n-1}\pm(Id_{C'}^{\otimes k}\otimes \Delta\otimes Id_{C'}^{\otimes l})\omega _{n}(e)\Big).$$
It follows from  property (3) of the family $\{\omega_{k}\mid k\geq 1\}$ that  $\widehat \omega d_{\cohoch}=d_{\cohoch}\widehat \omega$, since terms on either side of the equation arising from the internal differential match up in an obvious manner.  
\end{proof}

We can now define the type of highly structured coalgebras for which the coHochshild complex admits a natural comultiplication.

\begin{defn}\cite{hpst}  A \emph{weak Alexander-Whitney coalgebra} consists of a chain coalgebra $C$ such that the comultiplication $\Delta:C\to C\otimes C$ is a DCSH map, together with a choice of chain algebra map $\omega:\Om C\to \Om (C\otimes C)$ that realizes the DCSH structure of $\Delta$.  If the composite
$$\Om C \xrightarrow \omega \Om (C\otimes C) \xrightarrow q \Om C \otimes \Om C$$
is a coassociative comultiplication on $\Om C$, where $q$ denotes the Milgram equivalence, then $(C,\omega )$ is an \emph{Alexander-Whitney coalgebra}.  We call the composite $q\omega$ the \emph{associated loop comultiplication}.  
\end{defn}

If $\Delta:C\to C\otimes C$ is a DCSH map and $\omega:\Om C\to \Om (C\otimes C)$ realizes its DCSH structure, then $Id_{C}\otimes \Delta$ and $\Delta\otimes Id_{C}$ are both DCSH maps as well.  In particular, there are chain algebra maps 
$$Id_{C}\wedge\omega, \omega \wedge Id_{C}:\Om (C\otimes C)\to \Om (C\otimes C\otimes C)$$ realizing their DCSH structure, where the $k^{\text{th}}$-members of the associated families of $R$-linear maps, $(Id_{C}\wedge \omega)_{k}$ and $( \omega \wedge Id_{C})_{k}$, are given by following composites:
$$C\otimes C\xrightarrow{ \Delta ^{(k-1)}\otimes \omega_{k}} C^{\otimes k}\otimes (C\otimes C)^{\otimes k} \xrightarrow \cong (C\otimes C\otimes C)^{\otimes k}$$
and 

$$C\otimes C\xrightarrow{\omega_{k}\otimes \Delta ^{(k-1)}} (C\otimes C)^{\otimes k}\otimes C^{\otimes k} \xrightarrow \cong (C\otimes C\otimes C)^{\otimes k},$$
where the second map in each composite is the obvious permutation, and $\omega _{k}$ is the $k^{\text{th}}$-member of the family of $R$-linear maps associated to $\omega$.  For further justification of this construction, we refer the reader to section 1.1 in  \cite{hess2}.

\begin{defn}
A \emph{strict Alexander-Whitney coalgebra}  is a weak Alexander-Whitney coalgebra $(C,\omega)$ such that 
$$(Id_{C}\wedge\omega) \omega =(\omega\wedge Id_{C})\omega.$$
A \emph{quasistrict Alexander-Whitney coalgebra}  is a weak Alexander-Whitney coalgebra $(C,\omega)$ such that there is a derivation homotopy from  $(Id_{C}\wedge\omega) \omega$ to $(\omega\wedge Id_{C})\omega$.
\end{defn}

\begin{rmk}  Any strict Alexander-Whitney coalgebra is an Alexander-Whitney coalgebra, due to the naturality of the Milgram equivalence. 
\end{rmk}

\begin{ex} If $C$ is a cocommutative coalgebra, then $(C,\Om \Delta)$ is a strict Alexander-Whitney coalgebra.
\end{ex}

\begin{ex}\label{ex:hpst}  It was shown in \cite {hpst} that for  any $1$-reduced simplicial set $K$, there is a natural choice of chain algebra map $\omega _{K}:\Om C_{*}K\to \Om (C_{*}K\otimes C_{*}K)$ such that $(C_{*}K, \omega _{K})$ is an Alexander-Whitney coalgebra. We generalize this result in the Appendix, showing that if $K$  reduced, then $C_{*}K$ is a quasi-strict Alexander-Whitney coalgebra (Corollary \ref{cor:ck-quasistrict}). 

Let $K'$ be a simplicial set, and let $\E$ denote the simplicial suspension functor (cf., e.g., section 2.1(a) in \cite{hps2}).  Let $K=\E K'$.  Note that the generators of the free abelian group $C_{n+1}K$ are in natural, bijective correspondence with the generators of $C_{n}K'$, for all $n\geq 0$.  If $x$ is a generator of $C_{n}K'$, let $e(x)$ denote the corresponding generator of $C_{n+1}K$.

Let $\Delta$ denote the usual comultiplication on $C_{*}K'$.  Let $x\in C_{n}K'$, and write $\Delta (x)=x\otimes 1 + 1\otimes x + x_{i}\otimes x^{i}$, using Einstein summation notation.  It follows from the proof of Proposition 4.6 in \cite {hps2} that 
$$(\omega_{K})_{2}:C_{*}K\to (C_{*}K\otimes C_{*}K)^{\otimes 2}: e(x) \mapsto \pm \big(1\otimes e(x^{i})\big)\otimes\big(e(x_{i})\otimes 1\big), $$
where the sign follows from the Koszul rule, and that $(\omega_{K})_{n}=0$ for all $n\geq 3$.  Recall that $(\omega_{K})_{1}=\Delta$.  

In \cite{hps2} the authors concluded from Proposition 4.6 that if $K=\E K'$, then the associated loop comultiplication 
$$\psi_{K}:=q\omega_{K}:\Om C_{*}K\to \Om C_{*}K\otimes \Om C_{*}K$$ 
satisfies and is specified by
$$\psi_{K}\big(\si e(x)\big)=\si e(x)\otimes 1 + 1 \otimes \si e(x) +\si e(x_{i})\otimes \si e(x^{i}).$$

Note that in general $C_{*}K$ is not  a strict Alexander-Whitney coalgebra. 
\end{ex}

Having provided families of interesting examples of Alexander-Whitney coalgebras, we now prove that their coHochschild complexes admit natural comultiplicative structure.

\begin{thm}\label{thm:cohoch-comult}  Let $(C,\omega)$ be a weak Alexander-Whitney coalgebra,and let $\psi =q\omega: \Om C\to \Om C\otimes \Om C$.  The coHochschild complex on $C$, $\cohoch (C)$, admits a natural comultiplication $\widehat \psi$ such that 
$$\xymatrix{\Om (C)\ar[d] \ar[r]^\psi&\Om C\otimes \Om C\ar[d]\\
\cohoch (C) \ar[d]\ar[r]^(0.4){\widehat \psi}&\cohoch (C)\otimes\cohoch (C')\ar[d]\\
C\ar[r]^\Delta &C\otimes C}$$
commutes, where the vertical arrows are the natural inclusions and projections. 
  Moreover, the comultiplication on $\cohoch (C)$ is coassociative (respectively, coassociative up to chain homotopy) if $(C,\omega)$ is a strict (respectively, quasistrict) Alexander-Whitney coalgebra.
\end{thm}

\begin{proof}  By Theorem \ref{thm:ext-nat-cohoch}, there is a natural, induced chain map $$\widehat\omega:\cohoch (C)\to \cohoch (C\otimes C),$$
extending $\omega$.
Define the comultiplication $\widehat \psi$ on $\cohoch (C)$ to be the composite
$$\cohoch (C) \xrightarrow {\widehat \omega} \cohoch (C\otimes C)=\cohoch(C\otimes C, C\otimes C)\xrightarrow {\hat q} \cohoch (C)\otimes \cohoch (C),$$
where $\hat q$ is the extended Milgram equivalence of Theorem \ref{thm:milgram-cohoch}. It follows easily from the formulas in the proof of Theorem \ref{thm:ext-nat-cohoch} that $\widehat \psi$ is coassociative (respectively, coassociative up to chain homotopy) if $(Id_{C}\wedge\omega) \omega =(\omega\wedge Id_{C})\omega$ (respectively, if there is a derivation homotopy from  $(Id_{C}\wedge\omega) \omega$ to $(\omega\wedge Id_{C})\omega$).  
\end{proof}

 The following corollary is an immediate consequence of applying Theorem \ref {thm:cohoch-comult} to Example \ref{ex:hpst}.

\begin{cor}  If $K$ is a reduced simplicial set, then $\cohoch(C_{*}K)$ admits a natural comultiplication, which is coassociative up to chain homotopy.
\end{cor}

\begin{ex}\label{ex:comult-susp} In the case of a simplicial suspension $K=\E K'$, the formulas from the proof of Theorem \ref{thm:ext-nat-cohoch} reduce dramatically, enabling us to give a simple and explicit formula for $\widehat \psi$ on $\cohoch(C_{*}K)$.

Let $x\in C_{n}K'$, and, as in Example \ref{ex:hpst}, let $e(x)$ denote the corresponding generator of $C_{n+1}K$. Write $\Delta (x)=x\otimes 1 +1\otimes x + x_{i}\otimes x^{i}$.  Let $w\in \Om C_{*}K$, and write $\psi (w)=w_{j}\otimes w^{j}\in \Om C_{*}K\otimes \Om C_{*}K$.  It follows from the formulas in the proof of Theorem  \ref{thm:ext-nat-cohoch} that
\begin{align*}
\widehat \psi\big(e(x)\otimes w)=& \big(e(x)\otimes w_{j}\big) \otimes (1\otimes w^{j})+(1\otimes w_{j}) \otimes \big(e(x)\otimes w^{j}\big)\\
&\pm \big(1\otimes \si e(x_{i})\cdot w_{j}\big)\otimes \big(e(x^{i})\otimes w^{j}\big)\\
& \pm \big(e(x_{i})\otimes w_{j}\big)\otimes \big(1\otimes w^{j}\cdot \si e(x^{i})\big),
\end{align*}
where the signs follow from the Koszul rule. Note that $\widehat\psi $ is therefore strictly coassociative. 

If $K'$ is itself a simplicial suspension, then $\Delta(x)=x\otimes 1 + 1 \otimes x$ for all $x\in C_{>0}L$, so that
$$\widehat \psi\big(e(x)\otimes w\big)= \big(e(x)\otimes w_{j}\big) \otimes (1\otimes w^{j})+(1\otimes w_{j}) \otimes \big(e(x)\otimes w^{j}\big),$$
for all $w\in \Om C_{*}K$.  In other words, if $K$ is a double suspension, then the comultiplication on $\cohoch(C_{*}K)$ is the usual, unperturbed comultiplication on a tensor product of coalgebras.
\end{ex}

For our topological application in section \ref{sec:fls}, we need the following relative version of Theorem \ref{thm:cohoch-comult}.  Given two morphisms of coalgebras $f,g:C\to C'$, recall that ${}_{f}C_{g}$ denotes the $C'$-bicomodule structure on $C$ induced by $f$ on the left and $g$ on the right. 

\begin{thm}\label{thm:rel-cohoch-comult} Let $(C,\omega)$ and $(C',\omega')$ be weak Alexander-Whitney coalgebras, and let $\psi=q\omega$ and $\psi'=q\omega'$.  Let $f,g:C\to C'$ be morphisms of coalgebras that commute with the Alexander-Whitney structure, i.e., such that
$$\xymatrix{\Om C \ar [d]^\omega \ar [r]^{\Om f} &\Om C'\ar [d]^{\omega '}\\ \Om (C\otimes C)\ar[r]^{\Om (f\otimes f)}&\Om (C'\otimes C')}\quad \text{and}\quad\xymatrix{\Om C \ar [d]^\omega \ar [r]^{\Om g} &\Om C'\ar [d]^{\omega '}\\ \Om (C\otimes C)\ar[r]^{\Om (g\otimes g)}&\Om (C'\otimes C')} $$
commute.  Then $\cohoch ({}_{f}C_{g}, C')$ admits a natural comultiplication $\widehat \psi_{f,g}$ such that
$$\xymatrix{\Om C'\ar[d] \ar[rr]^{\psi'}&&\Om C'\otimes \Om C'\ar[d]\\
\cohoch ({}_{f}C_{g}, C') \ar[d]\ar[rr]^(0.4){\widehat \psi_{f,g}}&&\cohoch ({}_{f}C_{g}, C')\otimes\cohoch ({}_{f}C_{g}, C')\ar[d]\\
C\ar[rr]^\Delta &&C\otimes C}$$
commutes, where the vertical arrows are the natural inclusions and projections.   Moreover, the comultiplication $\widehat \psi_{f,g}$ is coassociative (respectively, coassociative up to chain homotopy) if $(C,\omega)$ and $(C',\omega ')$ are strict (respectively, quasistrict) Alexander-Whitney coalgebras.
 \end{thm}

\begin{proof}  We begin by defining a chain map, natural in $f$ and $g$,
$$\omega*\omega':\cohoch ({}_{f}C_{g}, C')\to \cohoch ({}_{f}C_{g}\otimes{}_{f}C_{g}, C'\otimes C').$$
We can then define the comultiplication $\widehat \psi_{f,g}$ on $\cohoch ({}_{f}C_{g}, C')$ to be given by the composite
$$\cohoch ({}_{f}C_{g}, C')\xrightarrow{\omega*\omega'} \cohoch ({}_{f}C_{g}\otimes{}_{f}C_{g}, C'\otimes C')\xrightarrow {\hat q}\cohoch ({}_{f}C_{g}, C')\otimes \cohoch ({}_{f}C_{g}, C'), $$
where $\hat q$ is the extended Milgram equivalence (Theorem \ref{thm:milgram-cohoch}).  

As in the proof of Theorem \ref{thm:ext-nat-cohoch}, we unfold the structure of $\omega$ and consider the associated family of $R$-linear maps
$$\{\omega_{k}:\overline C \to \overline {(C\otimes C)}^{\otimes k}\mid k\geq 1\}.$$ 
Using Sweedler-type notation and suppressing summation, we write
$$\omega _{k}(c)=(c^k_{1,1}\otimes c^k_{1,2})\otimes \cdots \otimes (c^k_{k,1}\otimes c^k_{{k,2}}).$$

For any $c\otimes w\in C\otimes \Om C'$, set
\begin{align*}
\omega*&\omega'(c\otimes w)\\
=&\sum _{k\geq 1\atop 1\leq i\leq k}\pm (c^k_{i,1}\otimes c^k_{i,2})\\
&\qquad\qquad\otimes \si \big(f(c^k_{i+1,1})\otimes f(c^k_{i+1,2})\big)|\cdots|\si\big(f(c^k_{k,1})\otimes f(c^k_{k,2})\big)\cdot \omega'(w)\\
&\qquad\qquad\qquad\qquad \cdot \si\big(g(c^k_{1,1})\otimes g(c^k_{1,2})\big)|\cdots|\si\big(g(c^k_{i-1,1})\otimes g(c^k_{i-1,2})\big),
\end{align*}
where the signs are determined by the Koszul rule. It is an exercise similar to the proof of Theorem \ref{thm:ext-nat-cohoch} to show that $\omega*\omega'$ is a chain map.  

As in the proof of Theorem \ref{thm:cohoch-comult}, it follows easily from the formulas in the proof of Theorem \ref{thm:ext-nat-cohoch} that $\widehat\psi_{f,g}$ is coassociative (respectively, coassociative up to chain homotopy) if $(C,\omega)$ and $(C',\omega')$ are strict (respective, quasistrict) Alexander-Whitney coalgebras.
\end{proof}

Again calling upon the work of Gugenheim and Munkholm, we know that for any pair of simplicial maps $g,h:K\to L$, where $L$ is reduced, the induced coalgebra maps $C_{*}g, C_{*}h:C_{*}K\to C_{*}L$ respect the natural Alexander-Whitney structure of their source and target. The corollary below is therefore an immediate consequence of  Theorem \ref{thm:rel-cohoch-comult}.

\begin{cor}  For any pair of simplicial maps $g,h:K\to L$, where $L$ is reduced, the coHochschild complex $\cohoch (C_{*}K, C_{*}L)$ admits a natural comultiplication that is coassociative up to chain homotopy, where $C_{*}K$ is considered as a $C_{*}L$-bicomodule via $C_{*}g$ on the left and $C_{*}h$ on the right.
\end{cor}

\begin{ex}\label{ex:comult-coince} In the case of a simplicial suspension $K=\E K'$, the formulas from the proof above again reduce, giving rise to a simple and explicit formula for $\widehat \psi_{g,h}$ on $\cohoch(C_{*}K,C_{*}L)$.

Let $x\in C_{n}K'$, and, as in Example \ref{ex:hpst}, let $e(x)$ denote the corresponding generator of $C_{n+1}K$. Write $\Delta (x)=x\otimes 1 +1\otimes x + x_{i}\otimes x^{i}$.  Let $w\in \Om C_{*}L$, and write $\psi (w)=w_{j}\otimes w^{j}\in \Om C_{*}L\otimes \Om C_{*}L$.  It follows from the formulas in the proof of Theorem  \ref{thm:ext-nat-cohoch} that
\begin{align*}
\widehat \psi_{g,h}\big(e(x)\otimes w\big)=& (e(x)\otimes w_{j}) \otimes (1\otimes w^{j})+(1\otimes w_{j}) \otimes (e(x)\otimes w^{j})\\
&\pm \Big(1\otimes \si f\big(e(x_{i})\big)\cdot w_{j}\Big)\otimes \big(e(x^{i})\otimes w^{j}\big)\\
& \pm \big(e(x_{i})\otimes w_{j}\big)\otimes \Big(1\otimes w^{j}\cdot \si g\big(e(x^{i})\big)\Big),
\end{align*}
where the signs follow from the Koszul rule. Note that $\widehat\psi_{g,h}$ is therefore strictly coassociative. 

If $K'$ is itself a simplicial suspension, then $\Delta(x)=x\otimes 1 + 1 \otimes x$ for all $x\in C_{>0}K'$, so that
$$\widehat \psi\big(e(x)\otimes w)= \big(e(x)\otimes w_{j}\big) \otimes (1\otimes w^{j})+(1\otimes w_{j}) \otimes \big(e(x)\otimes w^{j}\big),$$
for all $w\in \Om C_{*}L$.  In other words, if $K$ is a double suspension, then the comultiplication on $\cohoch(C_{*}K,C_{*}L)$ is the usual, unperturbed comultiplication on a tensor product of coalgebras.
\end{ex}

\section{CoHochschild complexes and homotopy coincidence spaces}\label{sec:fls}

Given two continuous maps $g,h :X\to Y$, their \emph{coincidence space}, which we denote $E_{g,h}$, is the equalizer of $g$ and $h$, i.e.,
$$E_{g,h}=\{x\in X\mid g(x)=h(x)\},$$
topologized as a subspace of $X$.  Another useful point of view is that $E_{g,h}$ is given by pulling back $(g,h): X\to Y\times Y$ over the diagonal map $\Delta:Y \to Y\times Y$, i.e., there is a pullback diagram
$$\xymatrix{ E_{g,h}\ar[d] \ar [r]&Y\ar [d]^\Delta\\ X\ar [r]^(0.4){(g,h)}&Y\times Y.}$$
If $Y=X$, then $E_{g,Id_{X}}=\operatorname{Fix}(g)$, the space of fixed points of $g$.  Of course, $E_{Id_{X},Id_{X}}=X$.

A homotopy-invariant version of the coincidence space of two maps $g,h:X\to Y$ is their \emph{homotopy coincidence space}, which we denote $\op E_{g,h}$ and for which one model is given by the pullback diagram
$$\xymatrix{ \op E_{g,h}\ar[d] \ar [r]&Y^{I}\ar [d]^{(ev_{0},ev_{1})}\\ X\ar [r]^(0.4){(g,h)}&Y\times Y,}$$
where $Y^{I}$ is the space of unbased paths on $Y$ and $ev_{t}$ is the map evaluating path at $t$.
In other words,  
$$\op E_{g,h}=\big\{(x,\ell)\in X\times Y^{I}\mid \ell (0)=g(x), \ell (1)=h(x)\big\}.$$
If $Y=X$, then $\op E_{g,Id_{X}}=\operatorname{Fix}^{ho}(g)$, the space of homotopy fixed points of $g$.  

Note that $\op E_{Id_{X},Id_{X}}=\operatorname{Fix}^{ho}(Id_{X})$ is the space of free loops on $X$, denoted $\op L X$. More generally, if $X=U\times V$, where $U$ and $V$ are submanifolds of a manifold $Y$ and $g, h:U\times V\to Y$ are given by projection onto the first and second coordinates, respectively, then $\op E_{g,h}$ is then exactly the space of open strings in $Y$ starting in $U$ and ending in $V$. 

Homotopy coincidence spaces play an important role in the study of geodesics on Riemannian manifolds.  Let $M$ be a closed, compact Riemannian manifold.  A slight generalization of a result of Gromoll and Meyer \cite{gromoll-meyer} states that if there is a field $\Bbbk$ such that the set of $\{\dim _{\Bbbk} H_{*}(\op L M;\Bbbk)\}$ of $\Bbbk$-Betti numbers of the free loop space on $M$ is unbounded, then $M$ admits an infinite number of distinct prime geodesics.  In \cite {grove-tanaka} Grove and Tanaka generalized Gromoll and Meyer's result, showing that if $\varphi$ is an isometry of finite order on $M$, then $M$ admits an infinite numbers of distinct, prime, $\varphi$-invariant geodesics  if there is a field $\Bbbk$ such that the set of $\{\dim _{\Bbbk} H_{*}(\operatorname{Fix}^{ho}(\varphi);\Bbbk)\}$  is unbounded.  

Homotopy coincidence spaces also show up in the theory of $p$-compact groups, where they give rise to \emph{$p$-compact groups of Lie type}.  Given an outer automorphism $\alpha$ of a $p$-compact group $BX$, one considers $BX(\alpha)=\operatorname{Fix}^{ho}(\alpha)$ \cite{moller}.

A model that facilitates computation of the comultiplicative structure of the homology of homotopy coincidence spaces should therefore have interesting applications in both homotopy theory and geometry, and perhaps in string topology as well. In this section we prove that the coHochschild complex provides just such a model, including the comultiplicative structure. More precisely, we prove the following theorem.

\begin{thm}\label{thm:model-htpy-coinc} If $g,h:K\to L$ are simplicial maps, where $L$ is a reduced simplicial set, then there is a quasi-isomorphism of chain complexes
$$\cohoch (C_{*}K,C_{*}L)\xrightarrow\simeq S_{*}\op E_{|g|,|h|}$$
that respects comultiplication up to chain homotopy, where $\cohoch (C_{*}K,C_{*}L)$ denotes the coHochschild complex, with $C_{*}K$ considered as a $C_{*}L$-bicomodule via $C_{*}g$ on the left and via $C_{*}h$ on the right.
\end{thm}

\begin{ex}  Let $n\geq 1$, and let $K=S^{2n+1}$, the simplicial sphere of dimension $2n+1$ with only one nondegenerate simplex of positive dimension, in dimension $2n+1$.  It is an easy exercice to show that the differential in $\cohoch(C_{*}S^{2n+1})$ is exactly zero.  Moreover, $S^{2n+1}$ is a double simplicial suspension, which implies (Exercise \ref{ex:comult-susp}) that the comultiplication on $\cohoch(C_{*}S^{2n+1})$ is the ordinary tensor product comultplication. Upon dualization, Theorem \ref{thm:model-htpy-coinc} therefore implies the well-known result that $H^{*}\op LS^{2n+1}$ is isomorphic as an algebra to $H^*(S^{2n+1})\otimes H^*(\Om S^{2n+1})$.
\end{ex}

In general, if $K$ is a simplicial double suspension, the fact that the comultiplication on $\cohoch(C_{*}K,C_{*}L)$ is simply the usual comultiplication on a tensor product of coalgebras (Example \ref{ex:comult-coince})  should make calculating the comultiplicative structure of $H_{*} \op E_{|g|,|h|}$ straightforward, once the homology has been calculated as graded abelian group.

To prove Theorem \ref{thm:model-htpy-coinc}, we begin by describing a  simplicial model for homotopy coincidence spaces, generalizing  the simplicial free loop space model developed in \cite {hess}.  We build this model as a twisted cartesian product, so we first recall the definition of this notion.  For any pair of simplicial maps $g,h$ with reduced codomain, we then apply homological perturbation theory to obtain a relatively small, homotopy-coassociative chain coalgebra, denoted $\op T(g,h)$, that is quasi-isomorphic to the normalized chain coalgebra on the simplicial homotopy coincidence construction on $g$ and $h$.   Finally, we prove the existence of a quasi-isomorphism from $\cohoch (C_{*}K,C_{*}L)$ to $\op T(g,h)$ that respects the comultiplications up to chain homotopy.

\subsection{A simplicial model for homotopy coincidence spaces}

We construct the simplicial model of a homotopy coincidence space as a twisted cartesian product, so we begin by recalling the necessary definitions from the theory of simplicial sets.

\begin{defn} Let $K$ be a simplicial set and $G$ a simplicial group, 
where the neutral element in any dimension is noted $e$. A 
degree $-1$ map of graded sets $\tau :K\to G$ is a \emph{twisting 
function} if 
\begin{align*}
\del _{0}\tau (x)&=\bigl(\tau (\del _{0}x)\bigr)^{-1}\tau (\del 
_{1}x)\\
\del _{i}\tau (x)&=\tau (\del 
_{i+1}x)\quad i>0\\
s_{i}\tau (x)&=\tau (s_{i+1}x)\quad i\geq 0\\
\tau (s_{0}x)&=e
\end{align*}
for all $x\in K$.
\end{defn}

\begin{rmk}\label{rmk:twisting-fcn} Let $K$ be a reduced simplicial set, and let $GK$ denote the Kan simplicial loop group on $K$ \cite{may}.  Let $\bar x\in (GK)_{n-1}$ denote a free group generator, corresponding to $x\in K_{n}$. There is a universal, canonical twisting function $\tau_{K}:K\to GK$, given by $\tau_{K}(x)=\bar x$.
\end{rmk}

Given a twisting function $\tau :K\to G$, where $G$ operates on the left on a simplicial set $L$, 
we can construct a \emph{twisted cartesian product} of $K$ 
and $L$, denoted $K\times _{\tau}L$, which is a simplicial set such 
that
$(K\underset{\tau}\times L)_{n}=K_{n}\times L_{n}$, with faces and 
degeneracies given by
\begin{align*} 
\del _{0}(x,y)&=(\del _{0}x,\tau (x)\cdot \del _{0}y)\\
\del _{i}(x,y)&=(\del _{i}x,\del 
_{i}y)\quad i>0\\
s_{i}(x,y)&=(s _{i}x,s_{i}y)\quad i\geq 0.
\end{align*}
If $L$ is a Kan complex, then the projection $K\underset{\tau}\times L\to K$ is 
a Kan fibration \cite {may}.

We can now generalize the construction of the simplicial free loop space model in \cite {hess} to model homotopy coincidence spaces.

\begin{defn} Let $g,h:K\to L$ be simplicial maps, where $L$ is a reduced simplicial set. The \emph{canonical homotopy coincidence construction} on $g$ and $h$, denoted $\op E_{g,h}$, is the twisted cartesian product 
$K\underset {\tau_{g,h}}\times GL,$
 where 
 $$\tau_{g,h} =(Gg\circ \tau_{K},Gh\circ\tau_{K}):K\to GL\times GL,$$ 
 and $GL\times GL$ acts on $GL$ by $(v,w)\cdot u=uvw^{-1}$. 
 \end{defn}
 
 Note that the construction of $\op E_{g,h}$ is clearly natural, which enables us to formulate the following definition
 
 \begin{defn} Let $\cat {Pair}$ denote the category where the objects are all pairs of simplicial maps $g,h:K\to L$, for all simplicial sets $K$ and all reduced simplicial sets $L$.  A morphism in $\cat {Pair}$ from a pair $g,h:K\to L$ to a pair $g',h':K'\to L'$ consists of a pair of simplicial morphisms $a:K\to K'$ and $b:L\to L'$ such that $f'a=bf$ and $g'a=bg$. 

The \emph{homotopy coincidence functor} $\op E:\cat {Pair}\to \cat {sSet}$  associates to any pair $(g,h)$ its homotopy coincidence space $\op E_{g,h}$ and to any morphism $(a,b):(g,h)\to (g',h')$ the simplicial map
$$\op E_{a,b}:\op E_{g,h}\to \op E_{g',h'}:(x, \bar y_{1}^{\pm 1}\cdots \bar y_{k}^{\pm 1})\mapsto \big(a(x), \overline{b( y_{1})}^{\pm 1}\cdots \overline{b( y_{k})}^{\pm 1}).$$
\end{defn}

We now show that  the simplicial homotopy coincidence construction does indeed provide simplicial models for topological homotopy coincidence spaces.
 
\begin{prop}  There is a commutative diagram of simplicial maps 
\begin{equation}\label{eqn:simpl-fls-model}
\xymatrix{GL\ar[d]^\simeq\ar[rr]^j&&\op E_{g,h}\ar [d]^\simeq\ar[r]^q& K\ar [d]^\simeq\\
\mathcal S_\bullet |GL|\simeq \mathcal S_\bullet \Om |L|\ar [rr]^(0.6){\mathcal S_\bullet i}&&\mathcal S_\bullet \op E_{|g|,|h|}\ar[r]^{\mathcal S_\bullet e}&\mathcal S_\bullet |K|,}
\end{equation} 
where $j$ and $q$ are the natural inclusion and projection, the vertical maps are weak equivalences of simplicial sets, and 
$$\Om |L|\xrightarrow i\op E_{|g|,|h|}\xrightarrow e |K|$$
is the obvious fibration sequence.  
\end{prop}

Consequently, $\H^*\op E_{g,h}$ and $\H^*\op E_{|g|,|h|}$ are isomorphic graded algebras.  For further details, we refer the reader to section 2.1 of \cite {hess}. The proofs there generalize easily  from the case of free loop spaces to the case of homotopy coincidence spaces. 

\subsection{Homological perturbation theory and a chain complex model for $\op E_{g,h}$}

The following classical notion from homological perturbation theory is necessary to our explanation of the relationship between $\op E_{|g|,|h|}$ and $\cohoch(C_{*}K,C_{*}L)$.

\begin{defn} Suppose that $\n :(X, \del ) \to  (Y,d)$ and
$f:(Y,d)\to (X,\del )$ are morphisms of (filtered) chain
complexes.  If $f\n = Id_{X}$  and there exists a (filtered)
chain homotopy
$h: (Y,d)\to  (Y,d)$ such that 
\begin{enumerate}
\item $dh +h d =\n f -Id_{Y}$,
\item $h \n =0$,
\item $fh =0$, and
\item $h \sp 2=0$,
\end{enumerate}
then $(X,d) \sdr{\n}f (Y,d)\circlearrowleft h$ is a (filtered) \emph{strong
deformation retract (SDR) of chain complexes.}  It is called \emph{Eilenberg-Zilber data} if $X$ and $Y$ are chain coalgebras and $\n$ is a morphism of coalgebras.
\end{defn}

\begin{rmk} If $(X,d) \sdr{\n}f (Y,d)\circlearrowleft h$ is Eilenberg-Zilber data, then 
$$(d\otimes Id_{X}+Id_{X}\otimes d)\bigl((f\otimes f)\Delta _{Y} h\bigr)+\bigl((f\otimes
f)\Delta_{Y}h\bigr)
d=\Delta_{X} f-(f\otimes f)\Delta_{Y},$$
i.e., $f$ is a map
of coalgebras up to chain homotopy. In fact, as stated precisely in the next theorem (due to Gugenheim and Munkholm and slightly strengthened in section 2.3 of \cite{hps2}), $f$ is a DCSH map, under reasonable local finiteness conditions.
\end{rmk}

\begin{thm}\label{thm:g-m} \cite{gugenheim-munkholm, hps2} Let $(X,d) \sdr{\n}f
(Y,d)\circlearrowleft h$ be Eilenberg-Zilber data such that $X$  and $Y$ are connected. Let $F_0=0$, and let $F_{1}$ be the composite
$$
\overline Y\xrightarrow f \overline X\xrightarrow{s^{-1}}s^{-1}\overline X.
$$
For $k\geq 2$, let  
$$
F_{k}=-\sum _{i+j=k}(F_{i}\otimes F_{j})\Delta _{Y} h:Y\rightarrow T^k(s^{-1}\overline X).
$$
If for all $y\in Y$, there exists $N(y)\in \mathbb N$ such that $F_{k}(y)=0$ for all $k>N(y)$, then
$$F=\prod _{k\geq 1}F_{k}=\bigoplus _{k\geq 1}F_{k}: Y \to \Om X$$
is a twisting cochain.
 In particular, $f:Y\to X$ is a DCSH map, and $\alpha_{F}:\Om Y\to \Om X$ realizes its strong homotopy structure.
 \end{thm}

\begin{ex}   Let
$K$ and $L$ be two simplicial sets.  There is  filtered Eilenberg-Zilber data
\begin{equation}\label{eqn:ez-sdr}
C_*K\otimes C_*L\sdr {\n_{K,L} }{f_{K,L}}C_*(K\times L)\circlearrowleft h_{K,L},
\end{equation}
where all structure is natural in $K$ and in $L$. The map $\nabla_{K,L} $ is the \emph{ shuffle (or Eilenberg-Zilber) map}, while
$f_{K,L}$ is the \emph{ Alexander-Whitney map}.  

If $K$ and $L$ are 1-reduced, then Theorem \ref{thm:g-m} implies immediately that $f_{K,L}$ is a DCSH map.  In section \ref{sec:EZ} we show, as a corollary of a more general result (Theorem \ref{thm:aw-dcsh}), that $f_{K,L}$ is in fact a DCSH map for all reduced simplicial sets $K$ and $L$. 
\end{ex}

The following result is the fundamental theorem of  homological perturbation theory.

\begin{thm}[The Basic Perturbation Lemma]\cite{r-brown}\label{thm:bpl} Let
$(X,\del) \sdr{\n}f (Y,d)\circlearrowleft h$ 
be a filtered SDR of chain complexes, where the filtrations are increasing and
bounded below. Let
$\theta:Y\to Y$ be a filtration-lowering linear map of degree $-1$  such that $(d+\theta)^2=0$. 
Define
\begin{align*}
 \n \sb {\infty}&=\n +\sum \sb {k>0}(h \theta)\sp k\n&f\sb {\infty}&=f+\sum\sb
{k>0}f(\theta h )\sp k\\
\\
\del \sb {\infty}&=\del +\sum\sb {k>0}f(\theta h )\sp
{k-1}\theta\n&h\sb {\infty}&=h +\sum \sb {k>0}h (\theta h
)\sp k \\
&&\hphantom{h\sb {\infty}}&=h +\sum \sb {k>0}(h \theta)\sp k h.
\end{align*}
Then $\del _{\infty}$, $\n _{\infty}$, $f_{\infty}$, and $h_{\infty}$ are all locally finite sums and 
$$(X,\del \sb {\infty}) \sdr{\n\sb {\infty}}{f\sb {\infty}}
(Y,d+\theta )\circlearrowleft h\sb {\infty}$$
is a filtered SDR with respect to the original filtrations of $X$ and $Y$.
\end{thm}

Hess proved an extended version of the Basic Perturbation Lemma in \cite [Theorem 4.1]{hess-ebpl}, explaining how to perturb a wide variety of algebraic structures in an SDR. The next theorem follows immediately from this Extended Basic Perturbation Lemma.

\begin{thm}\label{thm:perturb-EZdata}\cite{hess-ebpl} Let
$(X,\del) \sdr{\n}f (Y,d)\circlearrowleft h$ 
be filtered Eilenberg-Zilber data, where the filtrations are increasing and
bounded below. Let $\delta$ and $\Delta$ denote the comultiplications on $X$ and $Y$, respectively.  Given filtration-lowering linear maps $\theta: Y\to Y$ of degree $-1$ and $\zeta:Y\to Y\otimes Y$ of degree $0$ such that $(Y,\Delta +\zeta, d+\theta)$ is a chain coalgebra, there exists a chain map 
$\delta_{\infty}:(X,\del_{\infty})\to (X,\del_{\infty})\otimes (X,\del_{\infty})$ such that
\begin {enumerate}
\item $\delta_{\infty}-\delta$ is filtration-decreasing;
\item $\delta_{\infty}$ is coassociative up to chain homotopy; and
\item $f_{\infty}:(Y,\Delta +\zeta, d+\theta)\to (X,\delta_{\infty},\del _{\infty})$ is comultiplicative up to chain homotopy.
\end{enumerate}
\end{thm}

Explicit, natural formulas for $\delta _{\infty}$ and for the chain homotopy from $\delta_{\infty}f_{\infty}$ to $(f_{\infty}\otimes f_{\infty})(\Delta +\zeta)$ can be deduced from the formulas in Theorem 4.1 of \cite{hess-ebpl}. 

\begin{rmk}  It is probably true that $\delta_{\infty}$ is not only coassociative up to chain homotopy, but actually endows $X$ with the structure of an $A_{\infty}$-coalgebra.  Moreover, we expect that $\n_{\infty}$ and $f_{\infty}$ are $A_{\infty}$-coalgebra morphisms up to strong homotopy.  Since we did not require such powerful results in this article, we leave their proof to the ambitious reader.
\end{rmk}

As explained in Section 6.2 of \cite{hess-ebpl}, if $\tau:K\to G$ is a twisting function and $G$ acts on $L$, then the chain coalgebra $C_{*}(K\times _{\tau} L)$ is obtained by filtration-lowering perturbation of the differential and of the comultplication in $C_{*}(K\times L)$.  We can therefore apply Theorem \ref{thm:perturb-EZdata} to the Eilenberg-Zilber SDR (\ref{eqn:ez-sdr}) and establish the following result.

\begin{thm}\label{thm:brown-equiv} For each twisting function $\tau :K\to G$ and every simplicial set $L$ admitting a left action by $G$, 
there exists a homotopy-coassociative chain coalgebra $C_*K\widetilde\otimes_{\tau} C_*L$, obtained by filtration-lowering, natural perturbation of the differential and comultiplication of $C_*K\otimes C_*L$, together with an SDR
$$C_*K\widetilde\otimes_{\tau} C_*L \sdr{\n_{\tau}}{f_{\tau}}C_*(K\times 
_{\tau}L)\circlearrowleft\vp_{\tau},$$
where $\n_{\tau}$, $f_{\tau}$ and 
$\vp_{\tau}$ can be chosen naturally, and $\n _{\tau}$ and $f_{\tau}$ are comultiplicative up to natural chain homotopy.  
\end{thm}

Applying Theorem \ref{thm:brown-equiv} to the twisting function $\tau_{g,h}:K\to GL\times GL$ and to the conjugation action of $GL\times GL$ on $GL$, we obtain an SDR
$$C_*K\widetilde\otimes_{g,h} C_*GL \sdr{\n_{g,h}}{f_{g,h}}C_*(\op E_{g,h})\circlearrowleft\vp_{g,h},$$
where $\widetilde\otimes_{g,h}$, $\nabla _{g,h}$ , $f_{g,h}$ and $\vp_{g,h}$ are abbreviations for $\widetilde\otimes_{\tau_{g,h}}$, $\nabla_{\tau_{g,h}}$, $f_{\tau_{g,h}}$ and $\vp_{\tau_{g,h}}$.  Since all the constructions involved are natural in the pair $(g,h)$, there is a functor
\begin{equation}\label{eqn:twisted-model}
\op T:\cat {Pair} \to \cat {Coalg}_{R}^{hc},
\end{equation}
defined on objects by $\op T(g,h)=C_*K\widetilde\otimes_{g,h} C_*GL$, where $ \cat {Coalg}_{R}^{hc}$ is the category of homotopy-coassociative chain coalgebras and of homotopy-comultiplicative chain maps.

\subsection{The comparison map}

We now clarify the relationship between the two chain-level homotopy coincidence space models, $\op T(g,h)=C_*K\widetilde\otimes_{g,h} C_*GL$ and $\cohoch(C_{*}K,C_{*}L)$, by constructing a comparison map and showing that it is a quasi-isomorphism, respecting comultiplicative structure up to chain homotopy.

The foundation of our comparison map is the natural Szczarba equivalence
$$\alpha_{L}:\Om C_{*}L\xrightarrow \simeq C_{*}GL,$$
which is a quasi-isomorphism of chain algebras and a DCSH map for all reduced simplicial sets \cite{hpst}.

\begin{thm}\label{thm:cohoch-ez} Let $g,h:K\to L$ be simplicial maps, where $L$ is a reduced simplicial set.  Let  $\cohoch(C_{*}K,C_{*}L)$ denote the coHochschild complex, where $C_{*}K$ is considered as a $C_{*}L$-bicomodule via $C_{*}g$ on the left and $C_{*}h$ on the right.  Then there is a commutative diagram of chain complexes
\begin{equation}\label{eqn:cohoch-ez}
\xymatrix{
\Om C_{*}L\ar [d]_{\alpha_{L}}^\simeq\ar [r]^(0.4)\iota &\cohoch(C_{*}K,C_{*}L)\ar [d]_{\theta_{g,h}}^\simeq \ar @{->>}  [r]^(0.6){\pi_{\cohoch}} &C_{*}K\ar @{=}[d]\\
C_{*}GL\ar[r]^(0.4)\iota &C_*K\otimes _{{g,h}}C_*GL\ar @{->>}[r]^(0.6){\pi_{\op T}} &C_{*}K},
\end{equation}
where the vertical maps respect comultiplication, at least up to chain homotopy; $\theta_{g,h}$ is natural in the pair $(g,h)$; and the projection maps $\pi_{\cohoch}$ and $\pi_{\op T}$ are both defined in terms of the natural augmentation $C_{*}GL\to R$ sending the neutral element $e$ in degree $0$ to $1$ in $R$ and all other generators of $C_{0}GL$ to $0$.
\end{thm}

\begin{proof}   In Proposition 4.3 of \cite{hpst} the authors proved that $\alpha_{L}$ is  a morphism of chain coalgebras up to chain homotopy. The existence of $\theta_{g,h}$ and of its associated chain homotopy can be proved by a multi-stage inductive argument involving acyclic models.  

Define a partial order on the set 
$$(\mathbb Z_{+})^{\mathbb N}=\{(m_{1},...,m_{k})\mid k\geq 0, m_{i}\in \mathbb Z_{+} \; \forall i\}$$
of positive integer sequences (including the empty sequence) by
$$(m_{1},...,m_{k})< (n_{1},...,n_{l}) \Longleftrightarrow \begin{cases} \sum _{i=1}^k m_{i}<\sum _{j=1}^l n_{j}\\ \text{or}\\ \sum _{i=1}^k m_{i}=\sum _{j=1}^l n_{j}\quad\text{and}\quad k>l.\end{cases}$$
Given a graded $R$-module $V$ and $\vec m=(m_{1},...,m_{k})\in \mathbb Z_{+}^{\mathbb N}$, let $V_{\vec m}$ denote the $R$-module $V_{m_{1}-1}\otimes \cdots \otimes V_{m_{k}-1}$.
For any $m\in \mathbb Z_{+}$, let 
$$(\mathbb Z_{+})^{\mathbb N}_{m}=\{(m_{1},...,m_{k})\in (\mathbb Z_{+})^{\mathbb N}\mid \sum _{i=1}^km_{i}=m\}.$$
For any $m,k\in \mathbb Z_{+}$, let 
$$(\mathbb Z_{+})^{\mathbb N}_{m,k}=\{(m_{1},...,m_{j})\in (\mathbb Z_{+})^{\mathbb N}_{m}\mid j\geq k\}.$$
Note that 
$$(m_{1},...,m_{k})\in(\mathbb Z_{+})^{\mathbb N}_{m} \Longrightarrow k\leq m,$$
i.e., $(\mathbb Z_{+})^{\mathbb N}_{m,m+k}=(\mathbb Z_{+})^{\mathbb N}_{m,m}$ for all $k\geq 0$. On the other hand, $(\mathbb Z_{+})^{\mathbb N}_{m}=(\mathbb Z_{+})^{\mathbb N}_{m,1}$.

For any $n \in\mathbb Z_{+}$ and $\vec m\in (\mathbb Z_{+})^{\mathbb N}$, let $\cohoch (g,h)_{n; \vec m}$ denote the subcomplex
$$\Big(\big(C_{n}K\otimes \bigoplus_{\vec m'\leq \vec m}(\si C_{+}L)_{\vec m'}\big)\oplus \big(C_{<n}K\otimes T\si C_{+}L\big), d_{\cohoch}\Big)  $$
of the coHochschild complex $\cohoch (C_{*}K,C_{*}L)$.  Note that $\cohoch (g,h)_{n; \vec m}$ is also a sub coalgebra of $\cohoch (C_{*}K,C_{*}L)$.  Furthermore, given any simplicial maps $g,h:K\to L$, where $L$ is reduced,  there are comultiplicative isomorphisms
\begin{equation}\label{eqn:pullback-cohoch}
\cohoch (g,h)_{n;\vec m}\cong C_{*}K\underset {C_{*}(L\times L)}\square \cohoch (p_{1},p_{2})_{n;\vec m},
\end{equation}
while
\begin{equation}\label{eqn:pullback-twisted}
\op T(g,h)\cong C_{*}K\underset {C_{*}(L\times L)}\square\op T(p_{1},p_{2}),
\end{equation}
where $p_{1}, p_{2}:L\times L\to L$ denote the projections onto the first and second coordinate, respectively and 
$$M\underset {C}\square N$$ 
denotes the cotensor product of a  right $C$-comodule $M$ and a left $C$-comodule $N$ over a coalgebra $C$.

Let $\overline\Delta[n]$ denote the quotient of the standard simplicial $n$-simplex $\Delta[n]$ by its $0$-skeleton.   For all $\vec m=(m_{1},...,m_{k})\in(\mathbb Z_{+})^{\mathbb N}$,  let 
$$\overline \Delta [\vec m] =\overline\Delta[m_{1}]\vee\cdots\vee \overline\Delta[m_{k}].$$
For $n\geq 2$, let $\iota _{n}$ denote  the canonical generator of $C_{n}\Delta [n]$ and  let $\bar\iota _{n}$ denote both the canonical generator of $C_{n}\overline\Delta [n]$ and the image of this canonical generator in $C_{n}\overline \Delta [\vec m] $ when $m_{i}=n$ for some $i$. 

For any $n\in \mathbb Z_{+}$ and $\vec m\in(\mathbb Z_{+})^{\mathbb N}$, let 
$$j_{1}, j_{2}:\Delta [n]\to \overline\Delta [n]\vee\overline\Delta [\vec m ]\vee \overline\Delta [n]$$ 
denote the quotient map $\Delta [n]\to \overline \Delta [n]$ followed by the inclusion as the first, respectively last, summand of the wedge. Since $C_{*}G\overline\Delta[m]$ admits a contracting homotopy in positive degrees for all $m$, as proved by Morace and Prout\'e in \cite{morace-proutŽ}, an easy spectral sequence argument implies that
\begin{equation}
H_{k} \big(C_{*}\Delta[n]\otimes _{{j_{1}, j_{2}}}C_{*}G(\overline\Delta [n]\vee\overline\Delta [\vec m ]\vee \overline\Delta [n])\big) =0
\end{equation}
and therefore that 
\begin{equation}\label{eqn:acyclic}
H_{k} \big(\ker \pi_{\op T}) =0
\end{equation}
for all $k>0$.

Let $J\subset \mathbb N$, and let $V\subset (\mathbb Z_{+})^{\mathbb N}$.  For the purposes of this induction, we say that a collection of chain maps
$$\Theta_{J,V}=\{\theta_{g,h}^{n,\vec m}:\cohoch(g,h)_{n,\vec m}\to \op T(g,h)\mid n\in J, \vec m\in V, (g,h)\in \ob \cat{Pair}\} $$
is \emph{coherently admissible} if 
\begin{enumerate}
\item the collection is natural with respect to the pairs $(g,h)$;
\item the restrictions of $\theta_{g,h}^{n,\vec m}$ and of $\theta_{g,h}^{n',\vec m'}$ to the intersection of their domains are equal for all $n,n'\in J$ and all $\vec m, \vec m'\in V$;
\item the appropriate restrictions of diagram (\ref{eqn:cohoch-ez}) commute; 
\item the image of the restriction of $\theta _{g,h}^{n,\vec m}$ to $\ker \pi_{\cohoch}$ lies in $\ker \pi_{\op T}$ for all $(n,\vec m)\in J\times V$;
and
\item each $\theta_{g,h}^{n,\vec m}$ respects the comultiplications up to a natural chain homotopy
$$H_{g,h}^{n,\vec m}: \cohoch(g,h)_{n,\vec m}\to \op T(g,h),$$
 and the restrictions of $H_{g,h}^{n,\vec m}$ and of $H_{g,h}^{n'
 ,\vec m'}$ to the intersection of their domains are equal for all $n,n'\in J$ and all $\vec m, \vec m'\in V$.
\end{enumerate}

Note that the naturality of a coherently admissible collection $\Theta_{J,V}$ implies that for any simplicial maps $g,h:K\to L$, where $L$ is reduced, $\theta_{g,h}^{n,\vec m}$ is induced by the commuting diagram of left $C_{*}(L\times L)$-comodules
\begin{equation}\label{eqn:comod-diag}
\xymatrix{
\cohoch(p_{1},p_{2})_{n;\vec m}\ar [d]_{\theta_{p_{1},p_{2}}^{n,\vec m}} \ar @{->>}  [r]&C_{*}(L\times L)\ar @{=}[d]&&C_{*}K\ar[ll]_{C_{*}(g,h)}\ar@{=}[d]\\
\op T(p_{1},p_{2})\ar @{->>}[r]&C_{*}(L\times L)&&C_{*}K\ar[ll]_{C_{*}(g,h)}},
\end{equation}
for all $(n,\vec m )\in J\times V$, where $p_{1}, p_{2}:L\times L\to L$ are the projections onto the first and second factors, respectively. Here, we use the isomorphisms (\ref{eqn:pullback-cohoch}) and (\ref{eqn:pullback-twisted}).

We say furthermore that 
\begin{itemize}
\item condition $\mathbf{CA}_{n}$ is satisfied if there is a coherently admissible collection $\Theta_{J,V}$, where $J=\{ k\in \mathbb N\mid k\leq n\}$ and $V=(\mathbb Z_{+})^{\mathbb N}$;
\item condition $\mathbf{CA}_{n,m}$ is satisfied if there is a coherently admissible collection $\Theta_{J,V}$, where $J=\{ k\in \mathbb N\mid k\leq n\}$ and $V=(\mathbb Z_{+})^{\mathbb N}_{m}$;
\item condition $\mathbf{CA}_{n,m,l}$ is satisfied if there is a coherently admissible collection $\Theta_{J,V}$, where $J=\{ k\in \mathbb N\mid k\leq n\}$ and $V=(\mathbb Z_{+})^{\mathbb N}_{m,l}$.
\end{itemize}
Our goal is to prove that condition $\mathbf{CA}_{n}$ holds for all $n$, as this implies immediately that the desired chain map $\theta_{g,h}:\cohoch (C_{*}K,C_{*}L)\to \op T(g,h)$ exists, is natural in the pair $(g,h)$ and is comultiplicative up to chain homotopy.  

Assuming that condition $\mathbf{CA}_{n}$ holds for all $n$, it remains to show  that each $\theta_{g,h}$ is a quasi-isomorphism.  Consider first the case of the pair $(p_{1},p_{2})$. Composing $C_{*}p_{1}$ with $\pi_{\cohoch}$ and $\pi_{\op T}$, when $(g,h)=(p_{1},p_{2})$, we obtain a commutative diagram
\begin{equation}\label{eqn:cohoch-path}
\xymatrix{
\cohoch(C_{*}(L\times L),C_{*}L)\ar [d]_{\theta_{p_{1},p_{2}}} \ar @{->>}  [r]^(0.65)\simeq &C_{*}L\ar @{=}[d]\\
\op T(p_{1},p_{2})\ar @{->>}[r]^(0.65)\simeq &C_{*}L},
\end{equation}
in which the horizontal arrows are quasi-isomorphisms.  It follows that $\theta_{p_{1}, p_{2}}$ is a quasi-isomorphism as well. Since, as observed above, $\theta _{g,h}$ is induced by
\begin{equation}
\xymatrix{
\cohoch(p_{1},p_{2})\ar [d]_{\theta_{p_{1},p_{2}}} ^\simeq\ar @{->>}  [r]&C_{*}(L\times L)\ar @{=}[d]&&C_{*}K\ar[ll]_{C_{*}(g,h)}\ar@{=}[d]\\
\op T(p_{1},p_{2})\ar @{->>}[r]&C_{*}(L\times L)&&C_{*}K\ar[ll]_{C_{*}(g,h)}},
\end{equation} 
where both the source and the target of $\theta_{p_{1},p_{2}}$ are $C_{*}(L\times L)$-cofree, we conclude that every $\theta _{g,h}$ is a quasi-isomorphism.

To prove that condition $\mathbf{CA}_{n}$ holds for all $n$, we proceed inductively, establishing the following claims.
\begin{enumerate}
\item Condition $\mathbf{CA}_{0}$ holds such that $\theta_{g,h}\iota =\iota \alpha _{L}$ for all pairs $(g,h):K\to L\times L$ (cf., diagram (\ref{eqn:cohoch-ez})).
\item If condition $\mathbf{CA}_{n, m-1}$ holds, then condition $\mathbf{CA}_{n,m,m}$ holds.
\item If condition $\mathbf{CA}_{n, m,l+1}$ holds, then condition $\mathbf{CA}_{n,m,l}$ holds.
\item If condition $\mathbf{CA}_{n}$ holds, then condition $\mathbf{CA}_{n+1,0}$ holds.
\end{enumerate}
Since $(\mathbb Z_{+})^{\mathbb N}_{m}=(\mathbb Z_{+})^{\mathbb N}_{m,1}$ for all $m\in \mathbb Z_{+}$ and $(\mathbb Z_{+})^{\mathbb N}_{m,m+k}=(\mathbb Z_{+})^{\mathbb N}_{m,m}$ for all $k\geq 0$, it is clear that if this sequence of claims holds, then  condition $\mathbf{CA}_{n}$ holds for all $n$.
\medskip

\emph{Proof of Claim (1):}  It is easy to see that, since the differential on $\cohoch(C_{*}K, C_{*}L)_{0;\vec m}$ is untwisted for all $\vec m\in (\mathbb Z_{+})^{\mathbb N}$, we can choose
$\theta _{g,h}^{0,\vec m}$ to be the restriction of $Id_{C_{*}K}\otimes \alpha_{L}$ to $\cohoch(C_{*}K, C_{*}L)_{0;\vec m}$.
\medskip

\emph{Proof of Claim (2):}  Let $\vec m=(1,...,1)\in (\mathbb Z_{+})^{\mathbb N}_{m}$.
We now extend $\theta _{g,h}$ and $H_{g,h}$ naturally over $\cohoch(C_{*}K,C_{*}L)_{n,\vec m}$. Consider 
$$\iota _{n}\otimes \si \bar\iota _{{1}}|\cdots |\si\bar\iota_{{1}},$$
which is an element of 
$$\cohoch\big(j_{1},j_{2} )\big)_{n,\vec m},$$
where the bicomodule structure on $C_{*}\Delta[n]$ is induced by $C_{*}j_{1}$ on the left and by $C_{*}j_{2}$ on the right.
Note  that 
$d_{\cohoch}(\iota _{n}\otimes \si \bar\iota _{{1}}|\cdots |\si\bar\iota_{{1}})$ is an element of 
$$\Big(\sum _{\vec m'\in (\mathbb Z_{+})^{\mathbb N}_{m}\atop \vec m'<\vec m}\cohoch\big(C_{*}\Delta [n], C_{*}(\overline\Delta [n]\vee \overline \Delta [\vec m' ]\vee \overline \Delta [n] )\big)_{n,\vec m'}\Big)\bigcap\ker \pi_{\cohoch}.$$
Condition $\cat{CA}_{n,m-1}$ implies that $\theta_{j_{1}, j_{2}}\big(d_{\cohoch}(\iota _{n}\otimes \si \bar\iota _{{1}}|\cdots |\si\bar\iota_{{1}})\big)$ is defined and is a cycle in $\ker \pi_{\op T}$. If $n> 1$, it follows from (\ref{eqn:acyclic}) that there exists $\Phi\in  \ker \pi_{\op T}$ such that 
$$d_{\op T}\Phi=\theta_{j_{1}, j_{2}}d_{\cohoch}(\iota _{n}\otimes \si \bar\iota _{{1}}|\cdots |\si\bar\iota_{{1}}), $$
where $d_{\op T}$ denotes the differential in $\op T(g,h)$.
We can therefore set 
$$\theta_{j_{1}, j_{2}}(\iota _{n}\otimes \si \bar\iota _{{1}}|\cdots |\si\bar\iota_{{1}})=\Phi.$$

If $n=0$ or $n=1$, a simple calculation shows that we can set
$$\theta_{j_{1}, j_{2}}(\iota _{n}\otimes \si \bar\iota _{{1}}|\cdots |\si\bar\iota_{{1}})=\iota_{n}\otimes\alpha( \si \bar\iota _{{1}}|\cdots |\si\bar\iota_{{1}}).$$

For arbitrary $g, h:K\to L$ and generators $x\in C_{n}K$ and $y_{i}\in C_{1}L$ for $1\leq i\leq m$, 
let $\hat x:\Delta [n]\to K$, $\hat y_{i}: \overline \Delta [1]\to L$ denote the representing simplicial maps.  It is clear that $(\hat x, g\circ \hat x + \hat y_{1}+\cdots +\hat y_{k-1}+h\circ \hat x)$ is a morphism in $\cat {Pair}$ from $(j_{1},j_{2})$ to $(g,h)$, 
which implies that the pair $(\hat x, g\circ \hat x + \hat y_{1}+\cdots +\hat y_{k-1}+h\circ \hat x)$ induces chain maps
$$\cohoch(x; \vec y):\cohoch\big(C_{*}\Delta [n], C_{*}(\overline\Delta [n]\vee \overline \Delta [\vec m ]\vee \overline \Delta [n] )\big)\to \cohoch (C_{*}K, C_{*}L)$$
and
$$x*G\vec y:\op T(j_{1},j_{2}) \to \op T(g,h).$$
We set
$$\theta _{g,h}(x\otimes \si y_{1}|\cdots |\si y_{m})=(x*G\vec y)\circ \theta _{j_{1},j_{2}}(\iota _{n}\otimes \si \bar\iota _{{1}}|\cdots |\si\bar\iota_{1}),$$
so that $\theta_{g,h}\circ \cohoch (x;\vec y) =(x*G\vec y)\circ \theta _{j_{1},j_{2}}$, when applied to $\iota _{n}\otimes \si \bar\iota _{{1}}|\cdots |\si\bar\iota_{1}$.

Conditions (1)-(4) of the definition of a coherently admissible collection are then clearly satisfied, for $J=\{k\in \mathbb N\mid k\leq n\}$ and $V=(\mathbb Z_{+})^{\mathbb N}_{m,m}$.  Moreover, we can again call upon $(\ref{eqn:acyclic})$, in order to extend the collection of chain homotopies $H_{g,h}^{n;\vec m}$ naturally to $J=\{k\in \mathbb N\mid k\leq n\}$ and $V=(\mathbb Z_{+})^{\mathbb N}_{m,m}$ as well, thus fulfilling condition (5) of the definition of a coherently admissible collection.   In other words, condition $\mathbf{CA}_{n,m,m}$ holds.
\medskip 

\emph{Proof of Claim (3):}  The proof of this claim very closely resembles that of Claim (2).  We begin by using (\ref{eqn:acyclic}) to construct $\theta _{j_{1},j_{2}}^{n;\vec m}$ and its comultiplicativity chain homotopy $H_{j_{1},j_{2}}^{n;\vec m}$, for all $\vec m =(m_{1},...,m_{l})\in (\mathbb Z_{+})^{\mathbb N}_{m,l}$, where 
$$j_{1},j_{2}:\Delta[n]\to \overline\Delta [n]\vee \overline \Delta [\vec m ]\vee \overline \Delta [n]$$
are the usual quotient maps followed by inclusions.  We then extend to all pairs $(g,h)$ by naturality.
\medskip

\emph{Proof of Claim (4):} Consider the pair $j_{1},j_{2}:\Delta [n]\to \overline\Delta [n]\vee \overline \Delta [n]$. We need to show that if $\cat {CA}_{n}$ holds, then there exists $\Phi\in \ker( \pi_{\op T}:\op T(j_{1},j_{2})\to C_{*}\Delta [n])$ such that 
\begin{equation}\label{eqn:claim4}
d_{\op T}\Phi =\theta_{j_{1},j_{2}}d_{\cohoch}(\iota _{n}\otimes 1)-d_{\op T}(\iota_{n}\otimes e).
\end{equation}
We can then set $\theta_{g,h}^{n+1,\emptyset}(\iota _{n}\otimes 1)=\iota _{n}\otimes e + \Phi$, which implies that 
$$\pi_{\op T}\theta _{j_{1},j_{2}}^{n+1,\emptyset}(\iota _{n}\otimes 1)=\pi_{\cohoch}(\iota_{n}\otimes 1).$$  Extending by naturality, we obtain a collection of chain maps 
$$\{\theta_{g,h}^{n+1,\emptyset}: \cohoch (C_{*}K,C_{*}L)_{n+1; \emptyset}\to \op T(g,h)\mid (g,h)\in \ob \cat {Pair}\}$$
such that 
$$\xymatrix{
\cohoch(C_{*}K,C_{*}L)_{n+1,\emptyset}\ar [d]_{\theta_{g,h}} \ar @{->>}  [r]^(0.6){\pi_{\cohoch}} &C_{\leq n+1}K\ar [d]\\
\op T(g,h)\ar @{->>}[r]^(0.6){\pi_{\op T}} &C_{*}K}$$
commutes for all pairs $(g,h)$.
 
By condition $\cat {CA}_{n}$, the righthand side of equation (\ref{eqn:claim4}) is a well-defined cycle in $\ker \pi_{\op T}$.  Equation (\ref{eqn:acyclic}) therefore guarantees us the existence of the desired $\Phi$.  Again, a similar argument permit us to extend the comultiplicativity chain homotopy as well.

\medskip

Having now proved the four claims, we have completed the proof of the theorem. 

\end{proof}

The proof of Theorem \ref{thm:model-htpy-coinc} follows easily from the theorems above.

\begin{proof} [Proof of Theorem \ref{thm:model-htpy-coinc}] The desired quasi-isomorphism is equal to the following composite.
$$\xymatrix@1{\cohoch(C_{*}K,C_{*}L)\ar[r]_{\simeq}^{\theta_{g,h}}& C_*K\otimes _{t_{g,h}}C_*GL \ar[r]_(0.65){\simeq}^(0.65) {\nabla_{g,h}}& C_{*}\op E_{g,h}\ar[r]_{\simeq} &S_{*}\op E_{|g|,|h|}.}$$
Note that the first two maps in this composite are natural in $g$ and $h$.  The last map is obtained by applying the normalized chain functor to the middle vertical map in the diagram (\ref{eqn:simpl-fls-model}), which is not necessarily natural in $g$ and $h$.
\end{proof}

\section*{Appendix: The Eilenberg-MacLane SDR}\label{sec:EZ}
\setcounter{thm}{0}
\setcounter{equation}{0}
\renewcommand{\thethm}{A.\arabic{thm}}
\renewcommand{\thedefn}{A.\arabic{thm}}
\renewcommand{\theequation}{A.\arabic{equation}}

Our goal in this section is to prove a general existence result for DCSH maps, which implies in particular that the Alexander-Whitney map $$f_{K,L}:C_{*}(K\times L)\to C_{*}K\otimes C_{*}L$$ is a DCSH map, for all reduced simplicial sets $K$ and $L$.  It follows that $C_{*}K$ is naturally an Alexander-Whitney coalgebra for all reduced simplicial sets $K$.

We begin by an observation concerning the formula for the twisting cochain in Theorem \ref{thm:g-m} that proves useful in reaching our goal.

\begin{rmk} \label{rmk:Fk} Given Eilenberg-Zilber data $(X,d) \sdr{\n}f (Y,d)\circlearrowleft h$, there is a closed formula for each of the $F_{k}$'s in the statement of Theorem \ref{thm:g-m}.  For any $k\geq 2$, let 
$$h_{k}=\sum _{0\leq i\leq k-2} Id_{Y}^{\otimes i}\otimes \Delta_{Y}h\otimes Id_{Y}^{\otimes k-i-2}:Y^{\otimes k-1}\to Y^{\otimes k}$$
and let
\begin{equation}\label{eqn:Hk}
H_{k}=h_{k}\circ h_{k-1}\circ \cdots \circ h_{2}:Y\to Y^{\otimes k}.
\end{equation}
Then 
$$F_{k}=\sn {k+1} (\si f)^{\otimes k} \circ H_{k}.$$
\end{rmk}

In the development below, we use  the following helpful notation for simplicial expressions. 
 
 \begin{notn} If $J$ is any  set of non-negative integers
$j_1<j_2<\cdots <j_r$, let  $s_J$ denote the iterated degeneracy $s_{j_r}\cdots
s_{j_1}$, and let $| J|=r$.   

 For any $m\leq n\in \mathbb N$, let $[m,n]=\{j\in \mathbb N\mid m\leq j\leq n\}$. Let $\mathbf \Delta$ denote the category with objects
$$Ob \mathbf  \Delta=\{[0,n]\mid n\geq 0\}$$
and 
$$\mathbf  \Delta \bigl ([0,m], [0,n]\bigr)=\{ f:[0,m]\to [0,n]\;| \;f\text{ order-preserving set map}\}.$$
Viewing the simplicial $R$-module $M_{\bullet}$ as a contravariant functor from $\mathbf \Delta $ to the category of $R$-modules, given $x\in M_n:=M([0,n])$ and $0\leq a_1<a_2<\cdots <a_m\leq n$, let
$$x_{a_1...a_m}:=M(\mathbf a)(x)\in M_m$$
where $\mathbf a:[0,m]\to [0,n]:j\mapsto a_j$.
\end{notn}

Let $\mathcal A$ denote the usual functor from simplicial $R$-modules to $\cat {Ch}_{R}$, i.e., for any simplicial $R$-module $M_{\bullet}$,  the graded $R$-module underlying $\mathcal A(M_{\bullet})$ is $\{M_{n}\}_{n\geq 0}$, and the differential in degree $n$ is given by the alternating sum of the face maps from $M_{n}$ to $M_{n-1}$.  Let $\mathcal A_{N}$ denote its normalized variant. 

In Theorem 2.1a) of \cite{eilenberg-maclane} Eilenberg and MacLane gave explicit formulas for a natural SDR of chain complexes
\begin{equation}\label{eqn:em-sdr}
\mathcal A_{N}(M_{\bullet})\otimes \mathcal A_{N}(M'_{\bullet})\sdr \nabla f \mathcal A_{N}(M_{\bullet}\boxtimes M'_{\bullet})\circlearrowleft h,
\end{equation}
where $\boxtimes$ denotes the levelwise tensor product of simplicial $R$-modules.  In particular, if $x\in M_{m}$ and $x'\in M'_{n}$, then 
\begin{equation}\label{eqn:f}
f(x\boxtimes y)=\sum _{0\leq \ell \leq n} x_{0...\ell}\otimes y_{\ell...n}
\end{equation}
and
\begin{equation}\label{eqn:nabla}
\n (x\otimes x')=\sum _{0\leq \ell\leq n}\sum _{A\cup B=[0,n-1]\atop |A|=n-\ell, |B|=\ell}\pm s_{A}x\boxtimes s_{B}x',
\end{equation}
where the sign of a summand is the sign of the shuffle permutation corresponding to the pair $(A,B)$.

\begin{ex} If $R[K]$ denotes the free simplicial $R$-module generated by a simplicial set $K$, then $C_{*}K\otimes R\cong\mathcal A_{N}\big(R[K]\big)$. It follows that, when applied to $M_{\bullet}=R[K]$ and $M'_{\bullet}=R[L]$, for simplicial sets $K$ and $L$, Eilenberg and MacLane's strong deformation retract becomes the usual Eilenberg-Zilber/Alexander-Whitney equivalence
$$C_{*}K\otimes C_{*}L\sdr \nabla f  C_{*}(K\times L)\circlearrowleft h,$$
which is in fact Eilenberg-Zilber data. 
\end{ex}

\begin{rmk} Let $M_{\bullet}$ be a simplicial coalgebra over $R$, with levelwise comultiplication
$$\delta:M_{\bullet}\to M_{\bullet}\boxtimes M_{\bullet}.$$ 
The Eilenberg-MacLane SDR (\ref{eqn:em-sdr}) induces a coalgebra structure on $\op A_{N}(M_{\bullet})$, with coassociative comultiplication
$$\op A_{N}(M_{\bullet})\xrightarrow {\op A_{N}(\delta)} \op A_{N}(M_{\bullet}\boxtimes M_{\bullet})\xrightarrow f  \mathcal A_{N}(M_{\bullet})\otimes \mathcal A_{N}(M_{\bullet}).$$
Consequently, if $M_{\bullet}$ and $M_{\bullet}'$ are reduced simplicial coalgebras over $R$, then both $\op A_{N}(M_{\bullet }\boxtimes M'_{\bullet})$ and $\op A_{N} (M_{\bullet})\otimes  \op A_{N} (M'_{\bullet})$ are naturally chain coalgebras.
\end{rmk}

\begin{thm} \label{thm:aw-dcsh} If  $M_{\bullet}$ and $M_{\bullet}'$ are reduced simplicial coalgebras over $R$, with free underlying graded $R$-modules, then the Alexander-Whitney map 
$$f:\op A_{N}(M_{\bullet }\boxtimes M'_{\bullet}) \to \op A_{N} (M_{\bullet})\otimes  \op A_{N} (M'_{\bullet})$$
is a DCSH map.
\end{thm}

 To prove Theorem \ref{thm:aw-dcsh}, we apply Theorem \ref{thm:g-m} to the Eilenberg-MacLane SDR  (\ref{eqn:em-sdr}).  We must therefore prove local finiteness of the associated $F_{k}$'s, which follows from a technical result proved in \cite {hpst} (Lemma 5.3), expressed below in terms of simplicial $R$-modules instead of simplicial sets.

\begin{lem} \cite{hpst}\label{lem:hpst-htpy}  Let $M_{\bullet}$ and $M'_{\bullet}$ be simplicial $R$-modules.
Let $m<r\leq n$ be non-negative integers, and let $A$ and $B$ be disjoint sets of non-negative integers such that $ A\cup B=[m+1,n]$ and $|B|=r-m$.

Let $h^{A,B}: (M\boxtimes M')_{n}\to (M\boxtimes M')_{n+1}$ be the $R$-linear map given by
$$h^{A,B}(x\boxtimes x')=s_{A\cup\{m\}}\,x_{0\ldots r}\boxtimes s_B\,x'_{0\ldots mr\ldots n}.$$
for all $x\in M_{n}$ and $x'\in M'_{n}$.
Then the Eilenberg--MacLane homotopy in level $n$ 
$$h:\mathcal A_{n}(M_{\bullet}\boxtimes M'_{\bullet})=M_{n}\otimes M'_{n}\to M_{n+1}\otimes M'_{n+1}=\mathcal A_{n+1}(M_{\bullet}\boxtimes M'_{\bullet})$$ 
is given by
$$h(x\boxtimes x') =
\sum_{ m<r\\ A\cup B=[m+1,n] \atop
    |A|=n-r,\;|B|=r-m} 
\pm h^{A,B}(x\boxtimes x'),$$
where the sign corresponds to the sign of shuffle permutation associated to the couple $(A,B)$.
\end{lem}

\begin{proof}[Proof of Theorem \ref{thm:aw-dcsh}]  Fix $R$-bases $\mathsf B_{n}$ and $\mathsf B'_{n}$ of $M_{n}$ and $M'_{n}$, respectively, for all $n\geq 0$.  Define weight functions
$\zeta: \mathsf B_{n} \to \mathbb N  $ and $\zeta': \mathsf B_{n}' \to \mathbb N  $ as follows.  If $x\in \mathsf B_{n}$, then  
$$\zeta(x)=\max \{k\in \mathbb N\mid \exists (j_{1},...,j_{k})\in \mathbb N^k, y\in M_{n-k}\text { such that } x=s_{j_{k}}\cdots s_{j_{1}}y\}.$$
The function $\zeta'$ is defined similarly.

Let $x\in \mathsf B_{n}$ and $x'\in \mathsf B'_{n}$.  Observe that if $\zeta (x)+\zeta (x')>n$, then $x\boxtimes x'$ is necessarily a degenerate element of $M_{\bullet}\boxtimes M'_{\bullet}$.  In other words,
\begin{equation}\label{eqn:zeta-cond}
0\not = x\boxtimes x'\in \op A_{N}(M_{\bullet }\boxtimes M'_{\bullet})_{n}\Longrightarrow \zeta (x)+\zeta (x')\leq n.
\end{equation}

Consider the following bifiltration of $\op A_{N}(M_{\bullet }\boxtimes M'_{\bullet})$. For $p,n\geq 0$, let
$$\op F^{p,n}\big(\op A_{N}(M_{\bullet }\boxtimes M'_{\bullet})\big)$$ denote the graded submodule of $\op A_{N}(M_{\bullet }\boxtimes M'_{\bullet})$ that is generated by the set
$$\{x\boxtimes x'\mid x\in \mathsf B_{\leq n}, x'\in \mathsf B'_{\leq n}, \zeta (x)+\zeta (x')\geq p\}.$$
It follows from (\ref{eqn:zeta-cond}) that $\op F^{p,n}\big(\op A_{N}(M_{\bullet }\boxtimes M'_{\bullet})\big)=0$ for all $p>n$, so that we have a decreasing filtration
$$0\subset \op F^{n,n}\big(\op A_{N}(M_{\bullet }\boxtimes M'_{\bullet})\big)\subset \cdots \subset \op F^{1,n}\big(\op A_{N}(M_{\bullet }\boxtimes M'_{\bullet})\big)\subset \op A_{N}(M_{\bullet }\boxtimes M'_{\bullet}).$$ 

For any $k\geq 1$, consider the induced bifiltration 
{\smaller$$\op F^{p,n}\big(\op A_{N}(M_{\bullet }\boxtimes M'_{\bullet})^{\otimes k}\big)=\bigoplus _{p_{1}+\cdots +p_{k}=p\atop n_{1}+...+n_{k}=n} \mathcal F^{p_{1},n_{1}}\big (\op A_{N}(M_{\bullet }\boxtimes M'_{\bullet})\big)\otimes \cdots \otimes \mathcal F^{p_{k},n_{k}}\big (\op A_{N}(M_{\bullet }\boxtimes M'_{\bullet})\big).$$}
It is easy to check that the comultiplication 
$$\op A_{N}(M_{\bullet }\boxtimes M'_{\bullet})\to \op A_{N}(M_{\bullet }\boxtimes M'_{\bullet})\otimes \op A_{N}(M_{\bullet }\boxtimes M'_{\bullet})$$
is a bifiltered map. Moreover, it follows from implication (\ref{eqn:zeta-cond}) that
\begin{equation}\label{eqn:pn-cond}
\mathcal F^{p,n}\big (\op A_{N}(M_{\bullet }\boxtimes M'_{\bullet})^{\otimes k}\big)=0 \text{ for all }p>n\text { and } k\geq 1.
\end{equation}

To complete the proof of the theorem, we show that
\begin{equation}\label{eqn:zeta-h}
h\Big(\op F^{p,n}\big(\op A_{N}(M_{\bullet }\boxtimes M'_{\bullet})\big)\Big)\subset \op F^{p+2,n+1}\big(\op A_{N}(M_{\bullet }\boxtimes M'_{\bullet})\big)
\end{equation}
for all $p,n\geq 0$.
If (\ref{eqn:zeta-h}) holds, then
$$\Delta h\Big(\mathcal F^{p,n}\big (\op A_{N}(M_{\bullet }\boxtimes M'_{\bullet})\big)\Big)\subset \mathcal F^{p+2, n+1}\big (\op A_{N}(M_{\bullet }\boxtimes M'_{\bullet}) \otimes \op A_{N}(M_{\bullet }\boxtimes M'_{\bullet})\big),$$
which is the base step in an easy recursive argument showing that
$$H_{k+1}\Big(\mathcal F^{p,n}\big (\op A_{N}(M_{\bullet }\boxtimes M'_{\bullet})\Big)\subset \mathcal F^{p+2k, n+k}\big (\op A_{N}(M_{\bullet }\boxtimes M'_{\bullet})^{\otimes k}\big),$$
for all $k\geq 1$, where the map $H_{k}$ is defined as in (\ref{eqn:Hk}).

Equation (\ref{eqn:pn-cond}) therefore implies that for all $w\in \mathcal F^{p,n}\big(\op A_{N}(M_{\bullet }\boxtimes M'_{\bullet})\big)\big)$ and for all $k> n-p+1$,
$$F_{k}(w)=(\si f)^{\otimes k}\circ H_{k} (w)=0.$$
We have thus established the local finiteness of the $F_{k}$'s, which allows us to apply Theorem \ref{thm:g-m} and therefore conclude that $f:\op A_{N}(M_{\bullet }\boxtimes M'_{\bullet})\big)\to \op A_{N}(M_{\bullet })\otimes \op A_{n}(M'_{\bullet})$ is a DCSH map.

It remains only to verify (\ref{eqn:zeta-h}).  Let $x\in \mathsf B_{n}$ and  $x'\in \mathsf B'_{n}$.  If $\zeta (x)=l$ and $\zeta'(x')=m$, then there exist $y\in M_{n-l}$ and $y'\in M'_{n-m}$ such that $x_{0\cdots n}=y_{j_{0}\cdots j_{n}}$ and $x'_{0\cdots n}=y'_{j'_{0}\cdots j'_{n}}$, where $0\leq j_{0}\leq \cdots \leq j_{n}\leq n-l$ and $0\leq j'_{0}\leq \cdots \leq j'_{n}\leq n-m$.  It is clear that the level tensor product $x\boxtimes x'$ is degenerate if and only if there exists $k\in [0,n-1]$ such that $j_{k}=j_{k+1}$ and $j'_{k}=j'_{k+1}$.  On the other hand, if $l+m>n$, so that $(n-l)+(n-m)<n$, then the Pigeonhole Principle implies that  there exists $k\in [0,n-1]$ such that $j_{k}=j_{k+1}$ and $j'_{k}=j'_{k+1}$ and so $x\boxtimes x'$ is indeed degenerate.
\end{proof}

\begin{cor}\label{cor:ck-quasistrict}  If $K$ is a reduced simplicial set, then $C_{*}K$ is a quasi-strict Alexander-Whitney coalgebra.  If $K$ is $1$-reduced, then its associated loop comultiplication is strictly coassociative.
\end{cor}

\begin{proof} Theorem \ref{thm:aw-dcsh} implies that $C_{*}K$ is at least a weak Alexander-Whitney coalgebra, if $K$ is reduced.  To verify that it is quasi-strict, we apply an acyclic models argument.  

Let $\omega_{K}:\Om C_{*}K\to \Om (C_{*}K\otimes C_{*}K)$ denote the natural chain algebra map realizing the DCSH structure of the comultiplication on $C_{*}K$.  Both $(Id_{C_{*}K}\wedge\omega_{K}) \omega_{K}$ and $(\omega_{K}\wedge Id_{C_{*}K})\omega_{K}$ are natural in $K$ and therefore associated to families of $R$-linear morphisms
$$(\omega_{K}')_{k}, (\omega_{K}'')_{k}:C_{*}K\to (C_{*}K\otimes C_{*}K\otimes C_{*}K)^{\otimes k},$$
of degree $k-1$, for all $k\geq 1$, which are also natural in $K$.  Using as models the set $\{\overline \Delta[n]\mid n\geq 0\}$ of quotients of the standard $n$-simplices  by their $0$-skeletons, we can then prove by induction the existence of a family of $R$-linear morphisms
$$(\Phi_{K})_{k}:C_{*}K\to (C_{*}K\otimes C_{*}K\otimes C_{*}K)^{\otimes k},$$
of degree $k$, for all $k\geq 1$, and natural in $K$ giving rise to a derivation homotopy
$$\Phi=\sum _{k\geq 1}(s^{-1})^{\otimes k}\Phi _{k}:\Om C_{*}K\to \Om (C_{*}K\otimes C_{*}K\otimes C_{*}K)$$
from $(Id_{C_{*}K}\wedge\omega_{K}) \omega_{K}$ and $(\omega_{K}\wedge Id_{C_{*}K})\omega_{K}$. 

The $1$-reduced case was treated in \cite{hpst}.
\end{proof}

\begin{rmk}  It may be true that $C_{*}K$ is in fact an Alexander-Whitney coalgebra for all reduced $K$, i.e., that its associated loop comultiplication is strictly coassociative, but we have not yet proved it, as we did not need it for this article.
\end{rmk}

 \bibliographystyle{amsplain}
\bibliography{cohoch}
\end{document}